\documentclass [a4paper,12pt]{amsart}
\usepackage{amssymb}

\usepackage{graphicx}
\usepackage{epsfig}
\usepackage{amscd}

\newtheorem{theorem}{Theorem}
\theoremstyle{plain}

\newtheorem{definition}{Definition}
\newtheorem{example}{Example}

\newtheorem{lemma}{Lemma}

\newtheorem{proposition}{Proposition}
\newtheorem{remark}{Remark}

\numberwithin{equation}{section}
\numberwithin{theorem}{section}
\numberwithin{lemma}{section}
\numberwithin{proposition}{section}
\numberwithin{corollary}{section}
\numberwithin{remark}{section}
\numberwithin{example}{section}
\numberwithin{definition}{section}

\begin{document}
\title[ Local zeta functions \ and Newton polygons]{Local zeta function for curves, non degeneracy conditions and Newton polygons}
\author{M.J. Saia}
\address{Instituto de Matem\'{a}tica E Computa\c{c}ao, Universidade de S\~{a}o Paulo
at S\~{a}o Carlos, Av. do Trabalhador S\~{a}o -Carlense 400, CEP 13560-970,
S\~{a}o Carlos - SP, Brasil}
\email{mjsaia@icmc.sc.usp.br}
\thanks{The first named author was partially \ supported by CNPq-Grant 300556/92-6}
\author{W. A. Zuniga-Galindo}
\address{Department of Mathematics and Computer Science, Barry University, 11300 N.E.
Second Avenue\\
Miami Shores, Florida 33161, USA}
\email{wzuniga@mail.barry.edu}
\thanks{The second named author \ was supported by COLCIENCIAS-Grant \# 089-2000.
The second named author \ also thanks the partial support given by FAPESP
for visiting the Instituto de Matem\'{a}tica e Computa\c{c}ao, Universidade
de S\~{a}o Paulo, Campus S\~{a}o Carlos, \ in January 2000}

\begin{abstract}
This paper is dedicated \ to the description of the poles \ of the Igusa\
local zeta function $Z(s,f,v)$ when $f(x,y)$ satisfies \ a new non
degene\-racy condition called \textit{arithmetic non degeneracy}. More
precisely, we attach to each polynomial $f(x,y)$ a collection of convex sets 
$\Gamma ^{A}(f)=\left\{ \Gamma _{f,1},\cdots ,\Gamma _{f,l_{0}}\right\} $
called the \textit{arithmetic Newton polygon} of $f(x,y)$, and introduce the
notion of \ \textit{arithmetic\ non degeneracy with respect to }$\Gamma
^{A}(f)$. \ If $L_{v}$ is a $p-$adic field, and $f(x,y)\in L_{v}\left[ x,y%
\right] $ is arithmetically non degenerate, then the poles of $Z(s,f,v)$ can
be described explicitly in terms of the equations of the \ straight segments
that conform the boundaries of the convex sets $\Gamma _{f,1}$,$\cdots $, $%
\Gamma _{f,l_{0}}$. Moreover, \ the proof of main result gives an effective
procedure for computing $Z(s,f,v)$.
\end{abstract}

\maketitle

\section{Introduction}

Let $L_{v}$ be a $p-$adic field, and $v$\ the corresponding a
non-archimedean valuation. Let $O_{v}$ the ring of integers of \ $L_{v}$, $%
P_{v}$ the maximal ideal of $O_{v}$, and $\overline{L_{v}}=O_{v}/$ $P_{v}$\
the residue field of $L_{v}$. The cardinality of $\overline{L_{v}}$\ is
denoted by $q$, thus \ $\overline{L_{v}}=\mathbb{F}_{q}$. We fix a local
parameter $\pi $ for $O_{v}$. For $z\in L_{v}$, $\mid z\mid _{v}:=q^{-v(z)}$
denotes its normalized absolute value. Given a non-constant polynomial $%
g(x,y)\in L_{v}\mathbb{[}x,y]$, the Igusa local zeta function $Z(s,g,v)$
attached to $g(x,y)$ and $L_{v}$ is defined as

\begin{equation*}
Z(s,g,v):=\int\limits_{O_{v}^{2}}\mid g(x,y)\mid _{v}^{s}\mid dxdy\mid ,%
\text{ \ }
\end{equation*}
where $s\in \mathbb{C}$ satisfying $Re(s)>0$ and $\mid dxdy\mid $\ denotes \
a Haar measure of \ $L_{v\text{ }}^{2}$ normalized so that the measure of \ $%
O_{v}^{2}$\ \ is \ one.

Igusa showed that the local zeta function $Z(s,g,v)$ is a \ rational
function of $q^{-s}$, for general polynomials \ with coefficients in a $p-$%
adic field \cite{I1}. Igusa's local zeta functions are related to the number
\ of solutions of congruences modulo $p^{m}$ and to exponential sums modulo $%
p^{m}$. There are several conjectures and intriguing connections \ between
Igusa's local zeta functions with topology and singularity theory \ \cite{D1}%
.

This paper is dedicated \ to the description of the poles \ of the Igusa\
local zeta functions $Z(s,f,v)$ when $f(x,y)$ satisfies \ a new non
degeneracy condition called \textit{arithmetic non degeneracy}. More
precisely, we attach to each polynomial $f(x,y)$ a collection of convex sets 
$\Gamma ^{A}(f)=\left\{ \Gamma _{f,1},..,\Gamma _{f,l_{0}}\right\} $ called
the \textit{arithmetic Newton polygon} of $f(x,y)$, and introduce the notion
of \ \textit{arithmetic\ non degeneracy with respect to }$\Gamma ^{A}(f)$
(see section 4).

The set of degenerate polynomials, in the usual sense, is a proper subset of
the set of arithmetically non degenerate \ polynomials. In addition, for a
fixed geometric Newton polygon $\Gamma ^{geom}\subset \mathbb{R}^{2}$, the
set of degenerate complex polynomials in two variables with respect $\Gamma
^{geom}$ contains an open set, for the Zariski topology, consisting of
arithmetically non degenerate polynomials.

Our main result gives an explicit list for the possible poles of $Z(s,f,v)$,
when $f(x,y)\in L_{v}[x,y]$ is an arithmetically non degenerate polynomial,
in terms of the equations of the \ straight segments that conform the
boundaries of the convex sets $\Gamma _{f,1}$,.., $\Gamma _{f,l_{0}}$ (see
theorems \ref{theorem1}, \ref{main1}). Moreover, our proof gives an
effective procedure for computing $Z(s,f,v)$.

The notion of non degeneracy due Kouchnirenko \cite{Koc} is useful in the
computation of Milnor numbers, \ number of isolated solutions of polynomial
equations, etc. The authors hope that the \ notion of arithmetic non
degeneracy may be useful \ for other purposes, like those above mentioned,
different from the description of the poles of local zeta functions.

The \ description of the \ poles of the local zeta functions in terms of
Newton polygons is a well known fact for non degenerate polynomials in the
sense of Kouchnirenko \cite{Var}, \cite{LM}, \cite{D2}, \cite{D-H}, \cite
{Z-G1}, \cite{Z-G2}. In the degenerate case there are no previous results \
showing a relation between the poles of local zeta functions and Newton
polygons.

The poles of \ $Z(s,f,v)$, for a general polynomial $f$ \ in two variables,
can be described explicitly in terms of a resolution of singularities for $f$
\ (see \cite{V1}, \cite{V2}, \cite{M}, \cite{I3}, \cite{S}). The main result
of this paper is complementary to the mentioned results. Indeed, our main
result produces an explicit list of candidates for the poles of $Z(s,f,v)$
without using resolutions of singularities, but this result is only valid \
for a large class of polynomials.

In this paper, we work with $p-$adic fields of characteristic zero, however
all results are valid in positive characteristic.

This paper is organized as follows. In the section 2, we review some \
well-known results about geometric Newton polygons, non degeneracy
conditions in the sense of Kouchnirenko, and certain $p-$adic integrals
attached to non degenerate polynomials. In section 3, we \ compute
explicitly the local zeta function \ of\ some degenerate polynomials, in the
sense of Kouchnirenko. These examples constitute the basic models for the
effective computation of the local zeta functions of arithmetically non
degenerate polynomials. In section 4, we \ introduce the notion of
arithmetic Newton polygon and arithmetic non degeneracy condition. In
section 5, we compute \ the several $p-$adic integrals that appear in the
effective computation of the local zeta function of an arithmetically non
degenerate polynomial. In section 6, we prove the main result.

\section{Preliminaries}

\subsection{Geometric Newton polygons}

We set $\mathbb{R}_{+}=\{x\in \mathbb{R}\,\,\mid \,\,x\geqq 0\}$. Let $%
f(x,y)=\sum_{(i,j)}a_{i,j}x^{i}y^{j}$ be a non-constant polynomial in two
variables with coefficients in a field $K$ satisfying $f(0,0)=0$. The set $%
supp(f)=\left\{ (i,j)\in \mathbb{N}^{2}\,\mid \,a_{i,j}\neq 0\right\} $ is
called the \textit{support} of $f$. \textit{The geometric Newton polygon }$%
\Gamma ^{geom}(f)$ of $f$ is defined as the convex hull in $\mathbb{R}%
_{+}^{2}$ of the set $\cup _{(i,j)\in supp(f)}\left( (i,j)+\mathbb{R}%
_{+}^{2}\right) .$

{\ By a \textit{proper face} $\gamma $ of }$\Gamma ^{geom}(f)${, we mean a
non-empty convex set $\gamma $ obtained by intersecting }$\Gamma ^{geom}(f)$ 
{with an straight line $H$, such that }$\Gamma ^{geom}(f)$ {is contained in
one of two half-spaces determined by $H$. The straight line $H$ is called 
\textit{the supporting line} of $\gamma $. Thus the faces can be points
(zero-dimensional faces) or \ straight segments \ (one-dimensional faces).
The last ones are also called \textit{facets}. }

If $\gamma $ is a face of $\Gamma ^{geom}(f)$ {\ its \textit{face function}
is defined to be the polynomial } 
\begin{equation*}
f_{\gamma }(x,y):=\sum\limits_{(i,j)\in \gamma }a_{i,j}x^{i}y^{j}.
\end{equation*}

\begin{definition}
A non-constant polynomial $f(x,y)$ satisfying $f(0,0)=0$ is called\textit{\
non degenerate with respect to} $\Gamma ^{geom}(f)$, in the sense
Kouchnirenko, if it satisfies:

\begin{enumerate}
\item  the origin of \ $K^{2}$ is a singular point of \ $f(x,y)$;

\item  for each face $\gamma \subseteq $ $\Gamma ^{geom}(f)$, including $%
\Gamma ^{geom}(f)$ itself, the system of equations 
\begin{equation*}
f_{\gamma }(x,y)=\frac{\partial f_{\gamma }}{\partial x}(x,y)=\frac{\partial
f_{\gamma }}{\partial y}(x,y)=0,
\end{equation*}
does not have solutions on $K^{\times 2}$.
\end{enumerate}
\end{definition}

We set $\left\langle .,.\right\rangle $ for the usual inner product in $%
\mathbb{R}^{2}$, and identify the dual vector space with $\mathbb{R}^{2}$.
For $a\in \mathbb{R}_{+}^{2}$, we define

\begin{equation*}
m(a):=\inf_{x\in \Gamma (f)}\{\left\langle a,x\right\rangle \}.
\end{equation*}
Given any $a\in \mathbb{R}_{+}^{2}\smallsetminus \{0\}$ \textit{the first
meet locus }$F(a)$ of $a$ is defined as 
\begin{equation*}
F(a):=\{x\in \Gamma ^{geom}(f)\mid \,\left\langle a,x\right\rangle =m(a)\}.
\end{equation*}
The first meet locus $F(a)$ of $a\in \mathbb{R}_{+}^{2}\smallsetminus \{0\}$
is a proper face of $\Gamma ^{geom}(f)$.

Given a proper face $\gamma $ of $\Gamma ^{geom}(f)$, the \textit{cone
associated} to $\gamma $ is defined to be the set 
\begin{equation*}
\Delta _{\gamma }:=\left\{ a\in \mathbb{R}_{+}^{2}\smallsetminus \{0\}\,\mid
\,F(a)=\gamma \right\} .
\end{equation*}

We define an equivalence relation on $\mathbb{R}_{+}^{2}\smallsetminus \{0\}$
by 
\begin{equation*}
a\backsimeq a^{\prime }\,\,\,\text{if and only if}\,\,F(a)=F(a^{\prime }).
\end{equation*}

\noindent The \ equivalence classes of $\backsimeq $\ are the cones $\Delta
_{\gamma }$. The following two propositions describe the geometry of the
equivalences classes of $\backsimeq $ (see e.g. \cite{Var}, \cite{D2}).

\begin{proposition}
Let $\gamma $ be a proper face of $\Gamma ^{geom}(f)$. If $\gamma $ is a
one-dimensional face and \ $a$ is perpendicular vector on \ $\gamma $, then 
\begin{equation}
\Delta _{\gamma }=\{\alpha a\,\mid \,\alpha \in \mathbb{R},\,\ \alpha >0\}.
\label{cone1}
\end{equation}
If \ $\gamma $ is a zero-dimensional face and \ $w_{1},w_{2}$ are the facets
of $\Gamma ^{geom}(f)$ which contain $\gamma $, and $a_{1},a_{2}$ are
vectors which are perpendicular on respectively $w_{1},w_{2}$, then
\end{proposition}

\begin{equation}
\Delta _{\gamma }=\left\{\alpha _{1}a_{1}\,\,+\alpha _{2}a_{2}\,\mid
\,\alpha _{i}\in \mathbb{R},\,\ \alpha _{i}>0\right\}.  \label{cone2}
\end{equation}

A set of the form \ (\ref{cone2}) (respectively (\ref{cone1})) is called a 
\textit{\ cone strictly positive spanned} by the vectors $\left\{
a_{1},a_{2}\right\} $ (respectively by $\{a\}$). Let $\Delta $ be a strictly
positive spanned by a vector set $A$, with $A=\{a_{1},a_{2}\}$, or $A=\{a\}$%
. \ If $A$ is linearly independent over $\mathbb{R}$, the cone $\Delta $ is
called \textit{a simplicial cone}. In this last case, if \ the elements of $%
A $ are in $\mathbb{Z}^{2}$, the cone $\Delta $ is called {\ a \textit{%
rational simplicial }cone}. If $A$ can be completed to be a basis of $%
\mathbb{Z}-$module $\mathbb{Z}^{2}$, the cone $\Delta $ is called \textit{a
simple cone}.

A vector $a\in \mathbb{R}^{2}$ is called \textit{primitive} if the
components of $a$ are positive integers whose greatest common divisor is one.

For every facet of $\Gamma ^{geom}(f)$ there is a unique primitive vector in 
$\mathbb{R}^{2}$ which is perpendicular to this facet. We denote by $D$ the
set of such vectors. Thus each equivalence class under $\backsimeq $\ is a
rational simplicial cone spanned by elements of $D$.

From above considerations follow that there exists a partition of $\ \mathbb{%
R}_{+}^{2}$ \ of the form:

\begin{equation}
\mathbb{R}_{+}^{2}=\{(0,0)\}\bigcup \bigcup_{\gamma \subset \Gamma
^{geom}(f)}\Delta _{\gamma },  \label{cone3}
\end{equation}

\noindent where $\gamma $ runs through all proper faces of $\Gamma
^{geom}(f) $. \ We shall say that $\left\{ \Delta _{\gamma }\right\}
_{\gamma \subset \Gamma ^{geom}(f)}$ is a \textit{simplicial \ conical
subdivision }of $\mathbb{R}_{+}^{2}$\ subordinated to $\Gamma ^{geom}(f)$.

\subsection{\label{subsection} Local zeta functions and conical subdivisions}

Let $f(x,y)\in L_{v}[x,y]$ \ be a non-constant polynomial \ satisfying $%
f(0,0)=0$, and $\Gamma ^{geom}(f)$ its geometric Newton polygon. We fix \ a
simplicial \ conical subdivision $\left\{ \Delta _{\gamma }\right\} _{\gamma
\subset \Gamma ^{geom}(f)}$\textit{\ }of $\mathbb{R}_{+}^{2}$\ subordinated
to $\Gamma ^{geom}(f)$. We define 
\begin{equation*}
E_{\Delta _{\gamma }}:=\{(x,y)\in O_{v}^{2}\mid (v(x),v(y))\in \Delta
_{\gamma }\},
\end{equation*}
\begin{equation*}
Z(s,f,v,\Delta _{\gamma }):=\int\limits_{E_{\Delta _{\gamma }}}\mid
f(x,y)\mid _{v}^{s}\mid dxdy\mid ,
\end{equation*}
and 
\begin{equation*}
Z(s,f,v,O_{v}^{\times 2}):=\int\limits_{O_{v}^{\times 2}}\mid f(x,y)\mid
_{v}^{s}\mid dxdy\mid .
\end{equation*}
With the above notation, it follows from partition (\ref{cone3}) that 
\begin{equation}
Z(s,f,v)=Z(s,f,v,O_{v}^{\times 2})+\sum_{\gamma \subset \Gamma
^{geom}(f)}Z(s,f,v,\Delta _{\gamma }).  \label{cone5}
\end{equation}
Therefore the computation of the integral $\ Z(s,f,v)$ is reduced to the
computation of integrals of type $Z(s,f,v,O_{v}^{\times 2})$ \ and \ $%
Z(s,f,v,\Delta _{\gamma })$. In the case in which $f$ \ is non degenerate
with respect $\Gamma ^{geom}(f)$, the above mentioned integrals can be
computed explicitly using techniques of toroidal geometry, and the $p-$adic
stationary phase formula \cite[theorem 10.2.1]{I1}, \cite{D-H}, \cite{Z-G1}, 
\cite{Z-G2}, \cite{Z-G3}.

For a \ general $f$ \ formula \ (\ref{cone5}) is still valuable, if $\ f$
does not have singularities \ on $\left( L_{v}^{\times }\right) ^{2}$. By
applying \ the stationary phase formula (abbreviated SPF) repeatedly, it is
possible to compute effectively $Z(s,f,v,O_{v}^{\times 2})$ (see \ 
\cite[lemma 3.1]{Z-G1}, or \cite[lemma 2.4]{Z-G2}). For a future reference
we restate \ this result here.

\begin{lemma}
\label{lemma2}Let $g_{m}(x,y)=g_{0}(x,y)+\pi ^{m}g_{1}(x,y)\in L_{v}[x,y]$
be a \ non-constant polynomial, $m\geqslant 1$, such that $g_{m}(x,y)$,and $%
g_{0}(x,y)$ do not have singular points on $\left( L_{v}^{\times }\right)
^{2}$. There exists a constant $c(g_{0})$ such that if $m\geqq c(g_{0})$,
then 
\begin{equation*}
Z(s,g_{m},v,O_{v}^{\times 2})=\int\limits_{O_{v}^{\times 2}}\mid
g_{0}(x,y)\mid _{v}^{s}\mid dxdy\mid =\frac{U(q^{-s})}{1-q^{-1-s}},
\end{equation*}
with \ $U(q^{-s})\in \mathbb{Q[}q^{-s}\mathbb{]}$.
\end{lemma}

The computation of the integrals $Z(s,f,v,\Delta _{\gamma })$ presents two
cases. If $\ f_{\gamma }$ does not have singularities \ on $\left(
L_{v}^{\times }\right) ^{2}$, \ we shall say that $f$ \textit{is non
degenerate on }$\gamma $\textit{. }In this case the integrals $%
Z(s,f,v,\Delta _{\gamma })$ can be \ computed by using known techniques \cite
{LM}, \cite{D-H}, \cite{Z-G1}, \cite{Z-G2}. \ For a future reference we
restate this result here.

\begin{lemma}
\label{lemma3}Let $f(x,y)\in L_{v}[x,y]$ be a non-constant polynomial such
that $f(0,0)$ $=0$, and $f$ does not have singular points on $\left(
L_{v}^{\times }\right) ^{2}$.

\begin{enumerate}
\item  If $\Delta _{\gamma }$ is a two-dimensional cone, and $f_{\gamma
}(x,y)$ does not have singularities on the torus $\left( L_{v}^{\times
}\right) ^{2}$, then 
\begin{equation*}
Z(s,f,v,\Delta _{\gamma })=\text{ }U_{1}(q^{-s})\in \mathbb{Q[}q^{-s}\mathbb{%
]}.
\end{equation*}
\ 

\item  If $\Delta _{\gamma }$ is a one-dimensional cone, and $f_{\gamma
}(x,y)$ does not have singularities on $\left( L_{v}^{\times }\right) ^{2}$,
then 
\begin{equation*}
Z(s,f,v,\Delta _{\gamma })=\frac{U_{2}(q^{-s})}{(1-q^{-1-s})(1-q^{-(a+b)-d_{%
\gamma }s})}\text{, }U_{2}(q^{-s})\in \mathbb{Q[}q^{-s}\mathbb{]},
\end{equation*}

where $ax+by=d_{\gamma }$ is the equation of the supporting line of the
facet $\gamma $.
\end{enumerate}
\end{lemma}

The core of this paper is the explicit computation of \ integrals of\ type \ 
\begin{equation*}
Z(s,f,v,\Delta _{\gamma }),
\end{equation*}
when $f$ \ is \textit{degenerate\ on \ }$\Delta _{\gamma }$, but satisfies a
new non degeneracy condition.\ This will be done in sections 4, 5.

\begin{remark}
Given a $\Delta $ be a rational simplicial cone, there exists a partition of 
$\Delta $ into simple cones \cite[page 32-33]{K-K-M-S}. This fact can be
used to obtain a partition \ of $\mathbb{R}_{+}^{2}$\ \ in simple cones. In
this case, it is possible to compute the integral $Z(s,f,v,\Delta )$ \ in a
simple form, when $f$ is non degenerate, in the sense of Kouchnirenko.
However this procedure produces a larger list of candidates to poles for $%
Z(s,f,v,\Delta )$.
\end{remark}

\section{Model cases}

In this section we \ shall compute explicitly the local zeta functions for\
\ some degenerate polynomials, in the sense of Kouchnirenko. These examples
constitute the basic models for the effective computation of the local zeta
functions of arithmetically non degenerate polynomials.

\subsection{\label{sub1}The local zeta function of $%
(y^{3}-x^{2})^{2}+x^{4}y^{4}$}

We assume that the characteristic of the residue field of $L_{v}$ is
different from $2$. The origin of $L_{v}^{2}$ is the only singular \ point
of $f(x,y)=(y^{3}-x^{2})^{2}+x^{4}y^{4}$.This polynomial is degenerate with
respect to its geometric Newton polygon.

The conical subdivision of $\mathbb{R}_{+}^{2}$ subordinated to the
geometric Newton polygon of $f(x,y)$ is 
\begin{equation}
\mathbb{R}_{+}^{2}=\{(0,0)\}\bigcup \bigcup\limits_{j=1}^{9}\Delta _{j}\text{%
,}
\end{equation}
with

\begin{equation}
\begin{tabular}{|l|}
\hline
$\Delta _{1}:=(0,1)\mathbb{R}_{+}\text{ \ }$ \\ \hline
$\Delta _{2}:=(0,1)\mathbb{R}_{+}+(1,1)\mathbb{R}_{+}$ \\ \hline
$\Delta _{3}:=(1,1)\mathbb{R}_{+}$ \\ \hline
$\Delta _{4}:=(1,1)\mathbb{R}_{+}+(3,2)\mathbb{R}_{+}$ \\ \hline
$\Delta _{5}=:(3,2)\mathbb{R}_{+}\text{ }$ \\ \hline
$\Delta _{6}:=(3,2)\mathbb{R}_{+}+(2,1)\mathbb{R}_{+}$ \\ \hline
$\Delta _{7}:=(2,1)\mathbb{R}_{+}$ \\ \hline
$\Delta _{8}:=(2,1)\mathbb{R}_{+}+(1,0)\mathbb{R}_{+}$ \\ \hline
$\Delta _{9}:=(1,0)\mathbb{R}_{+}$ \\ \hline
\end{tabular}
\label{table0}
\end{equation}

In addition, all the \ $\Delta _{i}$ are rational simple cones. By using the
above subdivision is possible to reduce the computation of $Z(s,f,v)$ to the
computation of $Z(s,f,v,O_{v}^{\times 2}),$ and $Z(s,f,v,\Delta _{i}),$ $%
i=1,2,3,4,6,7,8,9$ (see\ref{subsection}).

\subsection{Computation of $Z(s,f,v,O_{v}^{\times 2})$}

Since $f(x,y)$ does not have singularities on $\left( L_{v}^{\times }\right)
^{2}$ \ the integral can be computed \ by using SPF \cite[theorem 10.2.1]{I1}%
:

\begin{equation*}
Z(s,f,v,O_{v}^{\times 2})=\left[ \left( 1-q^{-1}\right) ^{2}-q^{-2}N(f)%
\right] +\frac{q^{-2}(N(f)-1)\left( 1-q^{-1}\right) q^{-s}}{1-q^{-1-s}},
\end{equation*}
with $N(f)=Card\left( \left\{ (u,v)\in \left( \mathbb{F}_{q}^{\times
}\right) ^{2}\mid f(u,v)=0\right\} \right) $.

\subsection{Computation of $Z(s,f,v,\Delta _{i}),$ $i=1,2,3,4,6,7,8,9$}

These integrals \ correspond to the case in which $f$ is non degenerate on $%
\Delta _{i},$ $i=1$, $2$, $3$, $4$, $6$, $7$, $8$, $9$. \ The integrals
corresponding to $\Delta _{i},$ $i=1,2,$ \ can be calculated as follows.

\subsubsection{Case $Z(s,f,v,\Delta _{1})$}

\begin{eqnarray*}
Z(s,f,v,\Delta _{1}) &=&\sum_{n=1}^{\infty }\int\limits_{O_{v}^{\times
}\times \pi ^{n}O_{v}^{\times }}\mid f(x,y)\mid _{v}^{s}\mid dxdy\mid \\
&=&\sum_{n=1}^{\infty }q^{-n}\int\limits_{O_{v}^{\times 2}}\mid (\pi
^{3n}y^{3}-x^{2})^{2}+\pi ^{4n}x^{4}y^{4}\mid _{v}^{s}\mid dxdy\mid \\
&=&(1-q^{-1})^{2}\sum_{n=1}^{\infty }q^{-n}=q^{-1}(1-q^{-1}).
\end{eqnarray*}

\subsubsection{Case{\protect\Large \ } $Z(s,f,v,\Delta _{2})$}

\begin{eqnarray*}
Z(s,f,v,\Delta _{2}) &=&\sum_{m=1}^{\infty }\sum_{n=1}^{\infty
}\int\limits_{\pi ^{m}O_{v}^{\times }\times \pi ^{n+m}O_{v}^{\times }}\mid
f(x,y)\mid _{v}^{s}\mid dxdy\mid = \\
\sum_{m=1}^{\infty }\sum_{n=1}^{\infty
}q^{-2m-n-4ms}\int\limits_{O_{v}^{\times 2}} &\mid &(\pi
^{3n+m}y^{3}-x^{2})^{2}+\pi ^{4m+4n}x^{4}y^{4}\mid _{v}^{s}\mid dxdy\mid \\
&=&(1-q^{-1})q^{-1}\frac{q^{-2-4s}}{(1-q^{-2-4s})}.
\end{eqnarray*}

The other integrals can be computed in the same form. The following table
summarizes these computations:

\begin{equation}
\begin{tabular}{|l|l|}
\hline
Cone & $Z(s,f,v,\Delta )$ \\ \hline
$\Delta _{1}$ & $q^{-1}(1-q^{-1})$ \\ \hline
$\Delta _{2}$ & $(1-q^{-1})q^{-1}\frac{q^{-2-4s}}{(1-q^{-2-4s})}$ \\ \hline
$\Delta _{3}$ & $(1-q^{-1})^{2}\frac{q^{-2-4s}}{(1-q^{-2-4s})}$ \\ \hline
$\Delta _{4}$ & $(1-q^{-1})^{2}\frac{q^{-7-16s}}{(1-q^{-2-4s})(1-q^{-5-12s})}
$ \\ \hline
$\Delta _{6}$ & $\frac{(1-q^{-1})^{2}q^{-8-18s}}{(1-q^{-5-12s})(1-q^{-3-6s})}
$ \\ \hline
$\Delta _{7}$ & $\frac{(1-q^{-1})^{2}q^{-3-6s}}{(1-q^{-3-6s})}$ \\ \hline
$\Delta _{8}$ & $\frac{(1-q^{-1})q^{-4-6s}}{(1-q^{-3-6s})}$ \\ \hline
$\Delta _{9}$ & $q^{-1}(1-q^{-1})$ \\ \hline
\end{tabular}
\label{table}
\end{equation}

\subsection{Computation of $Z(s,f,v,\Delta _{5})$ (an integral on a
degenerate face in the sense of Kouchnirenko)}

\begin{eqnarray}
Z(s,f,v,\Delta _{5}) &=&\sum_{n=1}^{\infty }\int\limits_{\pi
^{3n}O_{v}^{\times }\times \pi ^{2n}O_{v}^{\times }}\mid f(x,y)\mid
_{v}^{s}\mid dxdy\mid  \notag \\
&=&\sum_{n=1}^{\infty }q^{-5n-12ns}\int\limits_{O_{v}^{\times 2}}\mid
(y^{3}-x^{2})^{2}+\pi ^{8n}x^{4}y^{4}\mid _{v}^{s}\mid dxdy\mid .
\label{2I7}
\end{eqnarray}

We set $f^{(n)}(x,y)=(y^{3}-x^{2})^{2}+\pi ^{8n}x^{4}y^{4}$, for $n\geqslant
1$. In order to compute (\ref{2I7}), we first compute the integral 
\begin{equation}
I(s,f^{(n)},v):=\int\limits_{O_{v}^{\times 2}}\mid (y^{3}-x^{2})^{2}+\pi
^{8n}x^{4}y^{4}\mid _{v}^{s}\mid dxdy\mid \text{, }n\geqslant 1.  \label{2I8}
\end{equation}
For that, we \ use the following change of variables 
\begin{equation}
\begin{array}{cccc}
\Phi : & O_{v}^{\times 2} & \longrightarrow & O_{v}^{\times 2} \\ 
& (x,y) & \longrightarrow & (x^{3}y,x^{2}y)\text{.}
\end{array}
\label{2I9}
\end{equation}

The map $\Phi $\ \ gives \ an analytic bijection of $O_{v}^{\times 2}$ onto
itself, and preserves the Haar measure because its Jacobian $J_{\Phi
}(x,y)=x^{4}y$ satisfies \ $\mid J_{\Phi }(x,y)\mid _{v}=1$, for every $%
x,y\in $ $O_{v}^{\times }$. Thus 
\begin{equation}
f^{(n)}\circ \Phi (x,y)=x^{12}y^{4}\widetilde{f^{(n)}}(x,y),  \label{2I10}
\end{equation}
with 
\begin{equation}
\widetilde{f^{(n)}}(x,y):=(y-1)^{2}+\pi ^{8n}x^{8}y^{4},  \label{2I11}
\end{equation}
and then from (\ref{2I9}), and (\ref{2I10}), it follows that

\begin{equation}
I(s,f^{(n)},v)=\int\limits_{O_{v}^{\times 2}}\mid \widetilde{f^{(n)}}%
(x,y)\mid _{v}^{s}\mid dxdy\mid .  \label{2I12}
\end{equation}

In order to \ compute the integral $I(s,f^{(n)},v)$, $n\geqslant 1$, we
decompose \ $O_{v}^{\times 2}$ as follows 
\begin{equation}
O_{v}^{\times 2}=\bigcup\limits_{y_{0}%
%TCIMACRO{
%\TeXButton{notequiv}{\not\equiv%
%}}%
%BeginExpansion
\not\equiv%
%
%EndExpansion
1\text{ mod }\pi }O_{v}^{\times }\times \{y_{0}+\pi O_{v}\}\bigcup
O_{v}^{\times }\times \{1+\pi O_{v}\},  \label{2I13}
\end{equation}
where $y_{0}$ runs through a set of representatives of $\mathbb{F}%
_{q}^{\times }$ in $O_{v}$. From \ partition \ (\ref{2I13}), and (\ref{2I11}%
), it follows that

\begin{equation}
I(s,f^{(n)},v)=(q-2)q^{-1}(1-q^{-1})+\int\limits_{O_{v}^{\times }\times
\{1+\pi O_{v}\}}\mid \widetilde{f^{(n)}}(x,y)\mid _{v}^{s}\mid dxdy\mid .
\label{2I14}
\end{equation}
The integral in the right side of (\ref{2I14}), admits the following
expansion:

\begin{equation}
I(s,f^{(n)},v)=(q-2)q^{-1}(1-q^{-1})+\sum_{j=1}^{\infty
}q^{-j}\int\limits_{O_{v}^{\times 2}}\mid \widetilde{f^{(n)}}(x,1+\pi
^{j}z)\mid _{v}^{s}\mid dxdz\mid ,  \label{2I15}
\end{equation}
\ where $\widetilde{f^{(n)}}(x,1+\pi ^{j}z)=\pi ^{2j}z^{2}+\pi
^{8n}x^{8}(1+\pi ^{j}z)^{4}$. In the next step, we compute \ explicitly $%
\mid \widetilde{f^{(n)}}(x,1+\pi ^{j}z)\mid _{v}^{s},$ \ for a fixed $n\geqq
1$. In this case, this can be done \ by a simple calculation:

\begin{equation}
\mid \widetilde{f^{(n)}}(x,1+\pi ^{j}z)\mid _{v}^{s}=\left\{ 
\begin{array}{c}
q^{-2js}\text{ \ \ \ \ if \ \ }1\leqq j\leqq 4n-1, \\ 
q^{-8ns}\mid z^{2}+x^{8}(1+\pi ^{j}z)^{4}\mid _{v}^{s}\text{ \ \ if \ \ }%
j=4n, \\ 
q^{-8ns}\text{ \ \ \ \ \ if \ \ }j\geqq 4n+1,
\end{array}
\right.  \label{2I17}
\end{equation}

for any $(x,z)\in O_{v}^{\times 2}$.

We note that the explicit computation of the absolute value of $\mid 
\widetilde{f^{(n)}}(x,1+\pi ^{j}z)\mid _{v}$ depends on an explicit
description \ of the set $\left\{ \left( w,z\right) \in \mathbb{R}^{2}\mid
w\leqq \min \{2z,8n\}\right\} ,$ for a fixed $n$.

Now from (\ref{2I15}), (\ref{2I17}) , and the identity $\sum%
\limits_{k=A}^{B}z^{k}=\frac{z^{A}-z^{B+1}}{1-z}$, we obtain an explicit
expression for $I(s,f_{n},v)$: 
\begin{eqnarray}
I(s,f^{(n)},v) &=&(q-2)q^{-1}(1-q^{-1})+(1-q^{-1})^{2}\frac{q^{-1-2s}}{%
(1-q^{-1-2s})}  \notag \\
&&-(1-q^{-1})^{2}\frac{q^{-4n-8ns}}{(1-q^{-1-2s})}+q^{-4n-8ns}I_{0}(s) 
\notag \\
&&+q^{-1}(1-q^{-1})q^{-4n-8ns},  \label{2I18}
\end{eqnarray}

where 
\begin{equation}
I_{0}(s):=\int\limits_{O_{v}^{\times 2}}\mid z^{2}+x^{8}(1+\pi
^{4n}z)^{4}\mid _{v}^{s}\mid dxdz\mid .  \label{2I19}
\end{equation}
Since the reduction modulo $\pi $ of \ the polynomial $\ z^{2}+x^{8}(1+\pi
^{4n}z)^{4}$ \ does not have singular points on $\mathbb{F}_{q}^{\times 2}$,
because the characteristic of the residue field is different of 2, SPF \
implies that $I_{0}(s)=\frac{U_{1}(q^{-s})}{(1-q^{-1-s})},$ where $%
U_{1}(q^{-s})$ \ is a polynomial in $q^{-s}$ independent of $j$ and $n$.

Finally, from (\ref{2I7}), (\ref{2I12}) and (\ref{2I18}), we obtain the
following explicit expression for $Z(s,f,v,\Delta _{5}):$

\begin{eqnarray}
Z(s,f,v,\Delta _{5}) &=&\frac{(q-2)(1-q^{-1})q^{-6-12s}}{(1-q^{-5-12s})}+%
\frac{(1-q^{-1})^{2}q^{-6-14s}}{(1-q^{-1-2s})(1-q^{-5-12s})}  \notag \\
&&-\frac{(1-q^{-1})^{2}q^{-9-20s}}{(1-q^{-1-2s})(1-q^{-9-20s})}+\frac{%
q^{-9-20s}U_{1}(q^{-s})}{(1-q^{-1-s})(1-q^{-9-20s})}  \notag \\
&&+\frac{(1-q^{-1})q^{-10-20s}}{(1-q^{-9-20s})}.  \label{2I20}
\end{eqnarray}

From the explicit formulas for $Z(s,f,v,\Delta _{i})$, $i=1,2,\cdots ,8,$\ \
and $Z(s,f,v,O_{v}^{\times 2})$ (see (\ref{table}), (\ref{2I20})), it
follows that the real parts of the poles of $Z(s,f,v)$ belong to the set

\begin{equation}
\{-1\}\bigcup \{-\frac{5}{12}\}\bigcup \{-\frac{1}{2},-\frac{9}{20}\}.
\label{poles1}
\end{equation}

\subsection{\label{sub2}The local zeta function of $%
(y^{3}-x^{2})^{2}(y^{3}-ax^{2})+x^{4}y^{4}$}

We set $g(x,y)=(y^{3}-x^{2})^{2}(y^{3}-ax^{2})+x^{4}y^{4}\in L_{v}[x,y]$,
with $a\in O_{v}^{\times }$, and $a%
%TCIMACRO{
%\TeXButton{notequiv}{\not\equiv%
%}}%
%BeginExpansion
\not\equiv%
%
%EndExpansion
1$ mod $\pi $. We assume that the characteristic of the residue field of $%
L_{v}$ is different from $2$. The polynomial $g$ is degenerate with respect
to its geometric Newton polygon. In addition, the origin of $L_{v}^{2}$ is
the only singular \ point of $g$.

By using the same decomposition of $\mathbb{R}_{+}^{2}$ into rational simple
cones as in the above example, it is possible to reduce the computation of $%
Z(s,g,v)$ to the computation of the $p-$adic integrals $Z(s,g,v,O_{v}^{%
\times 2}),$ \ $Z(s,g,v,\Delta _{i}),i=1$,$\cdots $, $8$, $9.$ \ 

\subsection{Computation of $Z(s,g,v,O_{v}^{\times 2})$}

Since $g(x,y)$ does not have singularities on $\left( L_{v}^{\times }\right)
^{2}$ \ the integral can be computed \ by using SPF \cite[theorem 10.2.1]{I1}%
:

\begin{equation*}
Z(s,g,v,O_{v}^{\times 2})=\left[ \left( 1-q^{-1}\right) ^{2}-q^{-2}N(g)%
\right] +\frac{q^{-2}(N(g)-1)\left( 1-q^{-1}\right) q^{-s}}{1-q^{-1-s}},
\end{equation*}
with $N(g)=Card\left( \left\{ (u,v)\in \left( \mathbb{F}_{q}^{\times
}\right) ^{2}\mid g(u,v)=0\right\} \right) $.

\subsection{Computation of $Z(s,g,v,\Delta _{i}),$ $i=1$, $2$, $3$, $4$, $6$%
, $7$, $8$, $9$.}

These integrals can be computed as the corresponding\ \ integrals in the
previous example. The following table summarizes these computations.

\begin{equation}
\begin{tabular}{|l|l|}
\hline
Cone & $Z(s,g,v,\Delta )$ \\ \hline
$\Delta _{1}$ & $q^{-1}(1-q^{-1})$ \\ \hline
$\Delta _{2}$ & $(1-q^{-1})\frac{q^{-3-6s}}{(1-q^{-2-6s})}$ \\ \hline
$\Delta _{3}$ & $(1-q^{-1})^{2}\frac{q^{-2-6s}}{(1-q^{-2-6s})}$ \\ \hline
$\Delta _{4}$ & $(1-q^{-1})^{2}\frac{q^{-7-24s}}{(1-q^{-2-6s})(1-q^{-5-18s})}
$ \\ \hline
$\Delta _{6}$ & $\frac{(1-q^{-1})^{2}q^{-8-27s}}{(1-q^{-5-18s})(1-q^{-3-9s})}
$ \\ \hline
$\Delta _{7}$ & $\frac{(1-q^{-1})^{2}q^{-3-9s}}{(1-q^{-3-9s})}$ \\ \hline
$\Delta _{8}$ & $\frac{(1-q^{-1})q^{-4-9s}}{(1-q^{-3-9s})}$ \\ \hline
$\Delta _{9}$ & $q^{-1}(1-q^{-1})$ \\ \hline
\end{tabular}
\label{table1}
\end{equation}

\subsection{Computation of $Z(s,f,v,\Delta _{5})$ (an integral on a
degenerate face in the sense of Kouchnirenko)}

\begin{eqnarray}
Z(s,g,v,\Delta _{5}) &=&\sum_{n=1}^{\infty }\int\limits_{\pi
^{3n}O_{v}^{\times }\times \pi ^{2n}O_{v}^{\times }}\mid g(x,y)\mid
_{v}^{s}\mid dxdy\mid  \notag \\
&=&\sum_{n=1}^{\infty }q^{-5n-18ns}\int\limits_{O_{v}^{\times 2}}\mid
(y^{3}-x^{2})^{2}(y^{3}-ax^{2})+\pi ^{2n}x^{4}y^{4}\mid _{v}^{s}\mid
dxdy\mid .  \label{2I21}
\end{eqnarray}

We set $g^{(n)}(x,y)=(y^{3}-x^{2})^{2}(y^{3}-ax^{2})+\pi ^{2n}x^{4}y^{4},$
for \ $n\geqslant 1$. In order to compute (\ref{2I21}), we first compute the
integral 
\begin{equation}
I(s,g^{(n)},v):=\int\limits_{O_{v}^{\times 2}}\mid
(y^{3}-x^{2})^{2}(y^{3}-ax^{2})+\pi ^{2n}x^{4}y^{4}\mid _{v}^{s}\mid
dxdy\mid \text{, }n\geqslant 1.  \label{2I22}
\end{equation}
By using the change of variables $\Phi (x,y)$ defined in (\ref{2I9}), we
have that

\begin{equation}
g^{(n)}\circ \Phi (x,y)=x^{18}y^{6}\widetilde{g^{(n)}}(x,y),  \label{2I23}
\end{equation}
with 
\begin{equation}
\widetilde{g^{(n)}}(x,y)=(y-1)^{2}(y-a)+\pi ^{2n}x^{2}y^{2}.  \label{2I24}
\end{equation}
Since $\Phi (x,y)$ gives a \ bijection of $O_{v}^{\times 2}$ onto itself
preserving the Haar measure, it follows from (\ref{2I22}), and (\ref{2I23}),
that

\begin{equation}
I(s,g^{(n)},v)=\int\limits_{O_{v}^{\times 2}}\mid \widetilde{g^{(n)}}%
(x,y)\mid _{v}^{s}\mid dxdy\mid .  \label{2I25}
\end{equation}

In order to compute the integral $I(s,g^{(n)},v)$, we decompose $%
O_{v}^{\times 2}$ \ as follows:

\begin{eqnarray}
O_{v}^{\times 2} &=&\left( O_{v}^{\times }\times \{y_{0}+\pi O_{v}\mid y_{0}%
%TCIMACRO{
%\TeXButton{notequiv}{\not\equiv%
%}}%
%BeginExpansion
\not\equiv%
%
%EndExpansion
1,a\text{ mod }\pi \}\right) \bigcup \left( O_{v}^{\times }\times \{1+\pi
O_{v}\}\right)  \notag \\
&&\bigcup \left( O_{v}^{\times }\times \{a+\pi O_{v}\}\right) .
\end{eqnarray}
From \ the above partition follows that

\begin{eqnarray}
I(s,g^{(n)},v) &=&(1-q^{-1})q^{-1}(q-3)  \notag \\
+\sum_{j=1\text{ }}^{\infty }q^{-j}\int\limits_{O_{v}^{\times }\times
O_{v}^{\times }} &\mid &\widetilde{g^{(n)}}(x,1+\pi ^{j}y)\mid _{v}^{s}\mid
dxdy\mid  \notag \\
+\sum_{j=1\text{ }}^{\infty }q^{-j}\int\limits_{O_{v}^{\times }\times
O_{v}^{\times }} &\mid &\widetilde{g^{(n)}}(x,a+\pi ^{j}y)\mid _{v}^{s}\mid
dxdy\mid ,  \label{2I27}
\end{eqnarray}
with 
\begin{equation}
\widetilde{g^{(n)}}(x,1+\pi ^{j}y)=\pi ^{2j}\gamma _{1,j}(x,y)y^{2}+\pi
^{2n}\gamma _{2,j}(x,y)x^{2},\text{ }  \label{2I28}
\end{equation}
\begin{equation}
\gamma _{1,j}(x,y)=1-a+\pi ^{j}y\text{, \ \ }\gamma _{2,j}(x,y)=(1+\pi
^{j}y)^{2},
\end{equation}
and 
\begin{equation}
\widetilde{g^{(n)}}(x,a+\pi ^{j}y)=\pi ^{j}\delta _{1}(x,y)y+\pi ^{2n}\delta
_{2}(x,y)x^{2},  \label{2I29}
\end{equation}
\begin{equation*}
\text{ }\delta _{1,j}(x,y)=(a-1+\pi ^{j}y)^{2}\text{, \ \ }\delta
_{2,j}(x,y)=(a+\pi ^{j}y)^{2}\text{.}
\end{equation*}
\ In the next step, we compute \ explicitly $\mid \widetilde{g^{(n)}}%
(x,1+\pi ^{j}y)\mid _{v}^{s}$, and \ $\mid \widetilde{g^{(n)}}(x,a+\pi
^{j}y)\mid _{v}^{s},$ for a fixed $n\geqq 1$. By some simple calculations,
it holds that

\begin{equation}
\mid \widetilde{g^{(n)}}(x,1+\pi ^{j}y)\mid _{v}^{s}=\left\{ 
\begin{array}{c}
q^{-2js}\text{ \ \ \ \ if \ \ }1\leqq j\leqq n-1, \\ 
q^{-2ns}\mid \gamma _{1,j}(x,y)y^{2}+\gamma _{2,j}(x,y)x^{2}\mid _{v}^{s}%
\text{ \ \ if \ \ }j=n, \\ 
q^{-2ns}\text{ \ \ \ \ \ if \ \ }j\geqq n+1,
\end{array}
\right.  \label{2I30}
\end{equation}

for any $(x,z)\in O_{v}^{\times 2}$, and 
\begin{equation}
\mid \widetilde{g^{(n)}}(x,a+\pi ^{j}y)\mid _{v}^{s}=\left\{ 
\begin{array}{c}
q^{-js}\text{ \ \ \ \ if \ \ }1\leqq j\leqq 2n-1, \\ 
q^{-2ns}\mid \delta _{1,j}(x,y)y+\delta _{2,j}(x,y)x^{2}\mid _{v}^{s}\text{
\ \ if \ \ }j=2n, \\ 
q^{-2ns}\text{ \ \ \ \ \ if \ \ }j\geqq 2n+1.
\end{array}
\right.  \label{2I31}
\end{equation}

We note that the explicit computation of the absolute values of $\mid 
\widetilde{g^{(n)}}(x,1+\pi ^{j}y)\mid _{v}$, respectively, $\mid \widetilde{%
g^{(n)}}(x,a+\pi ^{j}y)\mid _{v}$ depends on an explicit description \ of
the set 
\begin{equation*}
\left\{ \left( w,z\right) \in \mathbb{R}^{2}\mid w\leqq \min
\{2z,2n\}\right\} ,
\end{equation*}
respectively, of the set 
\begin{equation*}
\left\{ \left( w,z\right) \in \mathbb{R}^{2}\mid w\leqq \min
\{z,2n\}\right\} ,
\end{equation*}
for a fixed $n$.

Now, from \ \ (\ref{2I27}), (\ref{2I30}), (\ref{2I31}) and the identity $%
\sum\limits_{k=A}^{B}z^{k}=\frac{z^{A}-z^{B+1}}{1-z}$, we obtain the
following explicit expression for $I(s,g^{(n)},v)$:

\begin{eqnarray}
I(s,g^{(n)},v) &=&q^{-1}(1-q^{-1})(q-3)+\left( 1-q^{-1}\right) ^{2}\frac{%
q^{-1-2s}}{1-q^{-1-2s}}  \notag \\
&&-\left( 1-q^{-1}\right) ^{2}\frac{q^{-n-2ns}}{1-q^{-1-2s}}%
+q^{-n-2ns}I_{1}(s)+q^{-1}\left( 1-q^{-1}\right) q^{-n-2ns}  \notag \\
&&+\left( 1-q^{-1}\right) ^{2}\frac{q^{-1-s}}{1-q^{-1-s}}-\left(
1-q^{-1}\right) ^{2}\frac{q^{-2n-2ns}}{1-q^{-1-s}}+q^{-2n-2ns}I_{2}(s) 
\notag \\
&&+q^{-1}\left( 1-q^{-1}\right) q^{-2n-2ns}\text{,}  \label{2I32}
\end{eqnarray}

where 
\begin{eqnarray}
I_{1}(s) &:&=\int\limits_{O_{v}^{\times 2}}\mid \gamma
_{1,n}(x,y)y^{2}+\gamma _{2,n}(x,y)x^{2}\mid _{v}^{s}\mid dxdz\mid ,\text{ }
\notag \\
I_{2}(s) &:&=\int\limits_{O_{v}^{\times 2}}\mid \delta _{1,2n}(x,y)y+\delta
_{2,2n}(x,y)x^{2}\mid _{v}^{s}\mid dxdz\mid .\text{\ }  \label{2I33}
\end{eqnarray}
Since \ the reduction modulo $\pi $\ of the polynomials $\gamma
_{1,n}(x,y)y^{2}+\gamma _{2,n}(x,y)x^{2}$ and $\delta _{1,2n}(x,y)y+\delta
_{2,2n}(x,y)x^{2}$ does not have singular points on $\mathbb{F}_{q}^{\times
2}$, because the characteristic of the residue field is different from 2,
SPF \ implies that 
\begin{eqnarray}
I_{1}(s) &=&\frac{U_{1}(q^{-s})}{1-q^{-1-s}},\text{ }  \notag \\
I_{2}(s) &=&\frac{U_{2}(q^{-s})}{1-q^{-1-s}},\text{ }  \label{2I34}
\end{eqnarray}
where \ $U_{1}(q^{-s})$, and $U_{1}(q^{-s})$ are polynomials in $q^{-s}$,
independent of $j$, and $n$.

Finally, from (\ref{2I21}), and (\ref{2I32}), we obtain the following
explicit expression for $Z(s,f,v,\Delta _{5}):$

\begin{eqnarray}
Z(s,f,v,\Delta _{5}) &=&\frac{(q-3)(1-q^{-1})q^{-6-18s}}{(1-q^{-5-18s})}+%
\frac{(1-q^{-1})^{2}q^{-6-20s}}{(1-q^{-1-2s})(1-q^{-5-18s})}  \notag \\
&&-\frac{(1-q^{-1})^{2}q^{-6-20s}}{(1-q^{-1-2s})(1-q^{-6-20s})}+\frac{%
q^{-6-20s}U_{1}(q^{-s})}{(1-q^{-1-s})(1-q^{-6-20s})}  \notag \\
&&+\frac{(1-q^{-1})q^{-7-20s}}{(1-q^{-6-20s})}+\frac{(1-q^{-1})^{2}q^{-6-19s}%
}{(1-q^{-1-s})(1-q^{-5-18s})}  \notag \\
&&-\frac{(1-q^{-1})^{2}q^{-7-20s}}{(1-q^{-1-s})(1-q^{-7-20s})}+\frac{%
q^{-7-20s}U_{2}(q^{-s})}{(1-q^{-1-s})(1-q^{-7-20s})}  \notag \\
&&+\frac{(1-q^{-1})q^{-8-20s}}{(1-q^{-7-20s})}.  \label{2I35}
\end{eqnarray}

From the explicit formulas for $Z(s,f,v,\Delta _{i})$, $i=1,2,\cdots ,8,$\ \
and $Z(s,f,v,O_{v}^{\times 2})$ (see (\ref{table1}), (\ref{2I35})), it
follows that \ the real parts of the poles of $Z(s,f,v)$ belong to the set

\begin{equation}
\{-1\}\bigcup \{-\frac{5}{18},-\frac{1}{3}\}\bigcup \{-\frac{1}{2},-\frac{3}{%
10}\}\bigcup \{-\frac{7}{20}\}.  \label{poles2}
\end{equation}

\section{Arithmetic Newton polygons and \ arithmetic non degenerate
polynomials}

In this section we introduce the notions of arithmetic Newton polygon and
arithmetic non degeneracy for polynomials in two variables. The set of
arithmetically non degenerate polynomials\ contains strictly the set of non
degenerate polynomials in the sense of Kouchnirenko. In addition, an
explicit list of candidates for the poles of a local zeta function
associated to an arithmetically non degenerate polynomial \ can be given in
terms of the corresponding arithmetic Newton polygon (see section 6).

\subsection{ Arithmetically non degenerate semi-quasihomogeneous polynomials}

\begin{definition}
Let $K$ be a field, and $a$, $b$ two coprime positive integers. A polynomial 
$g(x,y)\in K[x,y]$ is quasihomogeneous with respect to the weight $(a,b)$,
if \ it has the form 
\begin{equation}
g(x,y)=cx^{u}y^{v}\prod_{i=1}^{l}\left( y^{a}-\alpha _{i}x^{b}\right)
^{e_{i}}\text{, \ }c\in K^{\times }\text{.}  \label{sem1}
\end{equation}
\end{definition}

From the definition follows that a quasihomogeneous polynomial $g(x,y)$
satisfies 
\begin{equation}
g(t^{a}x,t^{b}y)=t^{d}g(x,y)\text{, for every \ }t\in K^{\times }.
\label{sem2}
\end{equation}
The integer $d$ is called \ the \textit{weighted degree} of $g(x,y)$ with
respect to $(a,b)$. We note that our definition of quasihomogeneity \
coincides with the standard one (i.e. with (\ref{sem2})) after a finite
extension of $K$. \ 

\begin{definition}
\label{def4}A \ polynomial $g(x,y)\in K[x,y]$ is a semi-quasihomogeneous
polynomial with respect to the weight $(a,b)$, if \ it has the form 
\begin{equation}
g(x,y)=\sum\limits_{j=0}^{l_{g}}g_{j}(x,y),
\end{equation}
where each $g_{j}(x,y)$ is quasihomogeneous of weighted degree $d_{j}$ with
respect to $(a,b)$, and $d_{0}$ $<d_{1}<\cdots <$ $d_{l_{g}}$.\ 
\end{definition}

The polynomial $g_{0}(x,y)$\ is called \textit{the quasihomogeneous part }of
\ $g(x,y)$. If $\ $%
\begin{equation*}
g_{0}(x,y)=cx^{u_{0}}y^{v_{0}}\prod_{i=1}^{l_{0}}\left( y^{a}-\alpha
_{i,0}x^{b}\right) ^{e_{i,0}}
\end{equation*}
has singular points on the torus\ $\left( K^{\times }\right) ^{2}$, i.e. if $%
e_{i,0}>1$, for some $i=1,2,\cdots ,l_{0}$, then \ $g(x,y)$ is \ degenerate,
in the sense \ of Kouchnirenko, with respect to its geometric Newton polygon 
$\Gamma ^{geom}(g)$.

We put 
\begin{equation}
f_{j}(x,y):=c_{j}x^{u_{j}}y^{v_{j}}\prod_{i=1}^{l_{j}}\left( y^{a}-\alpha
_{i,j}x^{b}\right) ^{e_{i,j}}\text{, \ }c_{j}\in K^{\times }\text{,}
\label{kqs}
\end{equation}
$d_{j}$ for the weighted degree of \ $f_{j}(x,y)$ with respect to $(a,b)$,
i.e. 
\begin{equation*}
d_{j}:=ab\left( \sum\limits_{i=1}^{l_{j}}e_{i,j}\right) +au_{j}+bv_{j},
\end{equation*}
and \ 
\begin{equation}
f(x,y)=\sum\limits_{j=0}^{l_{f}}f_{j}(x,y),  \label{kqs1}
\end{equation}
with $d_{0}$ $<d_{1}<\cdots <$ $d_{l_{f}}$. From now on, we shall suppose
that $f(x,y)$\ is a degenerate polynomial in the sense of Kouchnirenko.

\begin{definition}
\label{4.1} Let $f(x,y)\in K[x,y]$ be a semi-quasihomogeneous polynomial, $\
\ \theta $ $\in K^{\times }$ a fixed root of $f_{0}(1,\theta ^{a})=0$, and $%
e_{j,\theta }$ the multiplicity of $\theta $\ as root of \ the polynomial
function \ $f_{j}(1,y^{a})$. To each $f_{j}(x,y)$ in (\ref{kqs1}), we
associate an straight line of the form 
\begin{equation*}
w_{j,\theta }(z):=(d_{j}-d_{0})+e_{j,\theta }z,\text{ }j=0,1,\cdots ,l_{f},
\end{equation*}
and define\textit{\ the arithmetic Newton polygon }$\Gamma _{f,\theta _{%
\text{ }}}$\textit{\ of} $f(x,y)$\ \textit{\ at } $\theta $\ as 
\begin{equation*}
\Gamma _{f,\theta _{\text{ }}}:=\{(z,w)\in \mathbb{R}_{+}^{2}\mid w\leqq
\min_{0\leq j\leq l_{f\text{ }}}\{w_{j,\theta }(z)\}.
\end{equation*}
The arithmetic Newton polygon $\Gamma ^{A}(f)$ of \ $f(x,y)$ is defined to
be the set 
\begin{equation*}
\Gamma ^{A}(f):=\{\Gamma _{f,\theta _{\text{ }}}\mid \theta \in K^{\times }%
\text{, \ }f_{0}(1,\theta ^{a})=0\}.
\end{equation*}
\end{definition}

We \ note that $\theta $ is one of the $\alpha _{i,0}$. In addition, $\Gamma
^{A}(f)$ \ is a finite collection of convex sets, and no \ just a convex set
as is $\Gamma ^{geom}(f)$.

If the polynomial $f_{0}(x,y)$ is interpreted as the weighted tangent cone
of $f(x,y)$, then for each tangent direction of the form \ $(\theta :1)$, $%
\theta \in K^{\times }$, there exists an arithmetic Newton polygon $\Gamma
_{f,\theta }$ at $\theta $, and all these Newton polygons constitute the
arithmetic Newton polygon $\Gamma ^{A}(f)$ of $f(x,y)$.

Let $\theta $ be a fixed root of $f_{0}(1,\theta ^{a})=0$. A point $Q$\ of
the topological boundary of $\ \Gamma _{f,\theta _{\text{ }}\mathit{\ }}$ is
called a \textit{vertex}, if \ $Q=(0,0)$ or if \ it is the intersection
point of two different straight lines of the form $w_{j,\theta }(z)$. Let $%
Q_{k}$, $k=0,1,\cdots ,r$, denote the vertices of the topological boundary
of $\Gamma _{f,\theta }$, with \ $Q_{0}:=(0,0)$. \ The boundary of $\Gamma
_{f,\theta }$ is conformed by $r$ straight segments, a half-line, and the
non-negative part of the horizontal axis of the $(w,z)-$plane. We \ put 
\begin{equation}
w_{k,\theta }(z)=(\mathcal{D}_{k}-d_{0})+\mathcal{E}_{k}z,\text{ \ }%
k=1,2,\cdots ,r,  \label{new2}
\end{equation}
for the equation of the straight segment between $Q_{k-1}$ and $Q_{k}$, and 
\begin{equation}
w_{r+1,\theta }(z)=(\mathcal{D}_{r+1}-d_{0})+\mathcal{E}_{r+1}z,
\end{equation}
for the equation of the half-line starting at $Q_{r}$. With this notation it
has that 
\begin{equation}
Q_{k}=(\tau _{k},(\mathcal{D}_{k}-d_{0})+\mathcal{E}_{k}\tau _{k}),\text{ }%
k=1,2,\cdots ,r,  \label{new3}
\end{equation}
with 
\begin{equation}
\tau _{k}:=\frac{\left( \mathcal{D}_{k+1}-\mathcal{D}_{k}\right) }{\mathcal{E%
}_{k}-\mathcal{E}_{k+1}}>0,\text{ }k=1,2,\cdots ,r.  \label{new4}
\end{equation}
We note that $\mathcal{D}_{k}=d_{j_{k}}$, and $\mathcal{E}%
_{k}=e_{j_{k},\theta }$, for some index $j_{k}$. In particular, $\mathcal{D}%
_{1}=d_{0}$, $\mathcal{E}_{1}=e_{0,\theta }$, and the first equation is $%
w_{1,\theta }(z)=\mathcal{E}_{1}z$. The slope \ $\mathcal{E}_{r+1}$ may be $%
0 $.

Let $Q$ \ be a vertex of the boundary of $\Gamma _{f,\theta _{\text{ }}}$. 
\textit{The face function }$f_{Q}(x,y)$ is defined to be the polynomial 
\begin{equation}
f_{Q}(x,y):=\sum\limits_{w_{j,\theta }(Q)=0}f_{j}(x,y),  \label{face}
\end{equation}
where $w_{j,\theta }(z)$\ is the straight line corresponding to $f_{j}(x,y)$.

Now we introduce a new notion of non degeneracy with respect to an
arithmetic Newton polygon.

\begin{definition}
A semi-quasihomogeneous polynomial $f(x,y)\in K\left[ x,y\right] $ \ is
called arithmetically non degenerate with respect to $\Gamma _{f,\theta _{%
\text{ }}}$at $\theta $, if it satisfies:

\begin{enumerate}
\item  the origin of $K^{2}$ is a singular point of $f(x,y)$;

\item  the polynomial $f(x,y)$ does not have singular points \ on $\left(
K^{\times }\right) ^{2}$;

\item  the system of equations 
\begin{equation*}
f_{Q}(x,y)=\frac{\partial f_{Q}}{\partial x}(x,y)=\frac{\partial f_{Q}}{%
\partial y}(x,y)=0
\end{equation*}
does not have solutions on $\left( K^{\times }\right) ^{2}$, for each vertex 
$Q$ $\neq (0,0)$ of \ boundary of $\Gamma _{f,\theta }$.
\end{enumerate}

A semi-quasihomogeneous polynomial $\ f(x,y)$\ $\in K\left[ x,y\right] $ is
called arithmetically non degenerate with respect to $\Gamma ^{A}(f)$, if \
it \ is arithmetically non degenerate with respect to $\Gamma _{f,\theta _{%
\text{ }}}$, \ for each $\theta \in K^{\times }$ satisfying $f_{0}(1,\theta
^{a})=0$.
\end{definition}

\subsection{ Arithmetically non degenerate polynomials}

Let $f(x,y)\in K[x,y]$ be a non-constant polynomial, and $\Gamma ^{geom}(f)$
its geometric Newton polygon. Each facet $\gamma \subset $ $\Gamma
^{geom}(f) $ has \ a perpendicular primitive vector $a_{\gamma }=$ $%
(a_{1}(\gamma ),a_{2}(\gamma ))$ and the corresponding supporting line has
an equation of the form $\left\langle a_{\gamma },x\right\rangle =d(\gamma )$%
.

\begin{definition}
Let $f(x,y)\in K[x,y]$ be a non-constant polynomial such that $f(x,y)$ is
semi-quasihomogeneous with respect to the weight $a_{\gamma }$, in the sense
of the definition \ref{def4}, for every\ facet $\gamma $ of $\ \Gamma
^{geom}(f)$. Let $\Gamma ^{A}(f,\gamma )$ be the \ arithmetic Newton polygon
of \ $f(x,y)$ considered as a semi-quasihomogeneous polynomial with respect
to the weight $a_{\gamma }$. The arithmetic Newton polygon $\Gamma ^{A}(f)$
of \ $f(x,y)$ is defined to be the set 
\begin{equation*}
\Gamma ^{A}(f):=\bigcup\limits_{\{\gamma \text{ a facet of }\Gamma
^{geom}(f)\}}\Gamma ^{A}(f,\gamma ).
\end{equation*}
\end{definition}

We note that the construction of the arithmetic polygon of $f$ needs that $%
f_{\gamma }$ can be factored as in (\ref{kqs}), for every facet $\gamma $ of 
$\Gamma ^{geom}(f).$ This condition can always be attained by passing to a
finite extension of $K$.

\begin{definition}
A polynomial $f(x,y)\in K[x,y]$ \ is called arithmetically non degenerate
with respect to its arithmetic Newton polygon, if \ for every facet $\gamma $
of $\ \Gamma ^{geom}(f)$, the semi-quasihomogeneous polynomial $f(x,y)$,
with respect to the weight $a_{\gamma }$, is arithmetically non degenerate
with respect to $\Gamma ^{A}(f,\gamma )$.
\end{definition}

We put $K=\mathbb{C}$, and fix a geometric Newton polygon $\Gamma
^{geom}\subset \mathbb{R}^{2}$. Then the set of degenerate polynomials in
two variables with respect $\Gamma ^{geom}$ contains an open set, for the
Zariski topology, consisting of arithmetically non degenerate polynomials.
The proof of this fact is similar to \cite[page 157]{A-GZ-V}.

\subsection{Examples}

The following three examples illustrate the above definitions for specific
polynomials, we assume that $K$ is a field of characteristic zero.

\begin{example}
\label{example1}We set $f(x,y)=(y^{3}-x^{2})^{2}+x^{4}y^{4}\in K[x,y]$. The
polynomial \ $f(x,y)$ is \ a degenerate, in the sense of Kouchnirenko,
semi-quasihomogeneous polynomial with respect to the weight $(3,2)$.The
origin of $K^{2}$ is the only singular point of $f(x,y)$. In this case the
arithmetic Newton polygon of $f(x,y)$ is equal to $\Gamma _{f,1}$. The
boundary of the Newton polygon $\Gamma _{f,1}$ is conformed by the straight
segments 
\begin{eqnarray*}
w_{1,1}(z) &=&2z,\text{ \ }0\leq z\leq 4, \\
w_{2,1}(z) &=&8,\text{ \ }z\geqslant 4,
\end{eqnarray*}
and the half-line \ $\left\{ \left( z,w\right) \in \mathbb{R}_{+}^{2}\mid
w=0\right\} $. Thus $\mathcal{D}_{1}=d_{0}=12$, $\mathcal{E}_{1}=2$,$%
\mathcal{D}_{2}=20$, $\tau _{0}=0$,$\ \tau _{1}=4$. The arithmetic \ Newton
polygon $\Gamma _{f,1}$ is showed in the next figure.

\begin{figure}[htp]
\hskip 4cm \epsfxsize=6cm \epsfbox{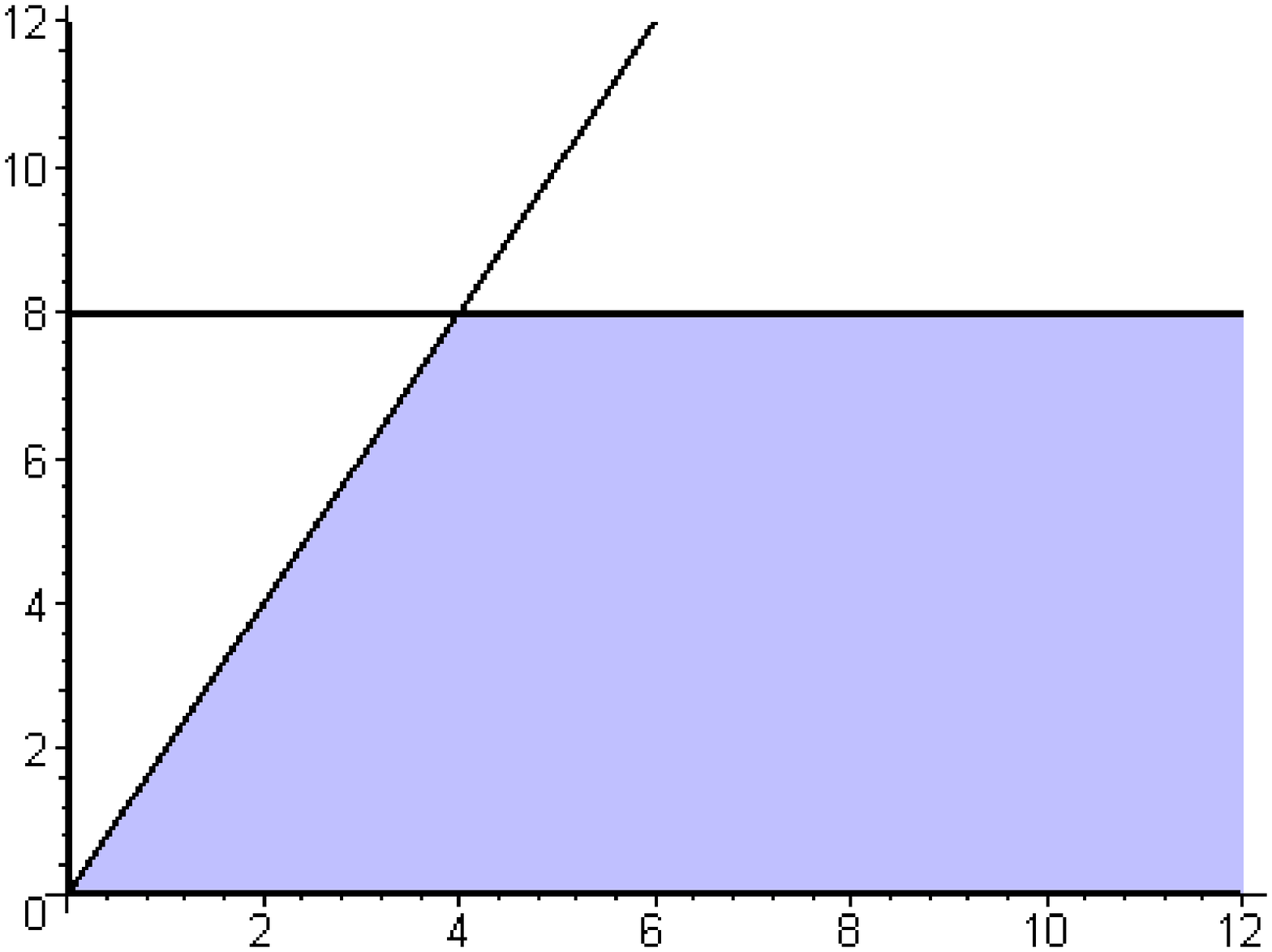}
\caption{$\Gamma ^{A}\left( f\right.)$}
\end{figure}

The face functions are 
\begin{equation*}
f_{(0,0)}(x,y)=(y^{3}-x^{2})^{2},\text{ }%
f_{(4,8)}(x,y)=(y^{3}-x^{2})^{2}+x^{4}y^{4}.
\end{equation*}
Since $f_{(4,8)}(x,y)$\ does not have singular points on $K^{\times 2}$, $%
f(x,y)$\ is arithmetically \ non degenerate.
\end{example}

\begin{example}
We set 
\begin{equation*}
g(x,y)=(y^{3}-x^{2})^{5}+(y^{3}-x^{2})^{3}x^{6}y^{3}+(y^{3}-x^{2})^{2}x^{12}+x^{24}\in K[x,y].
\end{equation*}
The polynomial \ $g(x,y)$ is \ a degenerate, in the sense of Kouchnirenko,\
semi-quasihomogeneous polynomial with respect to the weight $(3,2)$. The
origin of $K^{2}$ is the only singular point of $g(x,y)$. It can be
decomposed as 
\begin{equation*}
g(x,y)=g_{0}(x,y)+g_{1}(x,y)+g_{2}(x,y)+g_{3}(x,y),
\end{equation*}
where 
\begin{eqnarray*}
g_{0}(x,y) &=&(y^{3}-x^{2})^{5},\ \text{hence}\ d_{0}=30,\ e_{0,1}=5\ \text{%
and}\ w_{0,1}=5z; \\
g_{1}(x,y) &=&(y^{3}-x^{2})^{3}x^{6}y^{3},\ \text{hence}\ d_{1}=42,\
e_{1,1}=3\ \text{and}\ w_{1,1}=12+3z; \\
g_{2}(x,y) &=&(y^{3}-x^{2})^{2}x^{12},\ \text{hence}\ d_{2}=48,\ e_{2,1}=2\ 
\text{and}\ w_{2,1}=18+2z; \\
g_{3}(x,y) &=&x^{24},\ \text{hence}\ d_{3}=72,\ e_{3,1}=0\ \text{and}\
w_{3,1}=42.
\end{eqnarray*}
\textrm{\ }

Since $g_{0}(1,\theta ^{a})=0$ has only one root, $\theta =1$, the
arithmetic Newton polygon $\Gamma ^{A}\left( g\right) $ of $g(x,y)$ is equal
to $\Gamma _{g,1}$. The boundary of $\Gamma ^{A}\left( g\right) $ is
conformed by the straight segments 
\begin{eqnarray*}
w_{1,1}(z) &=&5z,\text{ \ }0\leq z\leq 6, \\
w_{2,1}(z) &=&18+2z,\text{ \ }6\leq z\leq 12, \\
w_{3,1}(z) &=&42,\text{ \ }z\geqslant 12,
\end{eqnarray*}

and the half-line \ $\left\{ \left( z,w\right) \in \mathbb{R}_{+}^{2}\mid
w=0\right\} $. Thus $\mathcal{D}_{1}=d_{0}=30$, $\mathcal{E}_{1}=5$, $%
\mathcal{D}_{2}=48$, $\mathcal{E}_{2}=2$, $\mathcal{D}_{3}=72$, $\mathcal{E}%
_{3}=0$, $\tau _{0}=0$,$\ \tau _{1}=6$, $\tau _{2}=12$.

We observe that the line $w(z)=12+3z$ meets $\Gamma _{g,1}$ only at $(6,30)$%
. The arithmetic Newton polygon of $g$ is showed in the next figure.

\begin{figure}[htp]
\hskip 4cm \epsfxsize=6cm \epsfbox{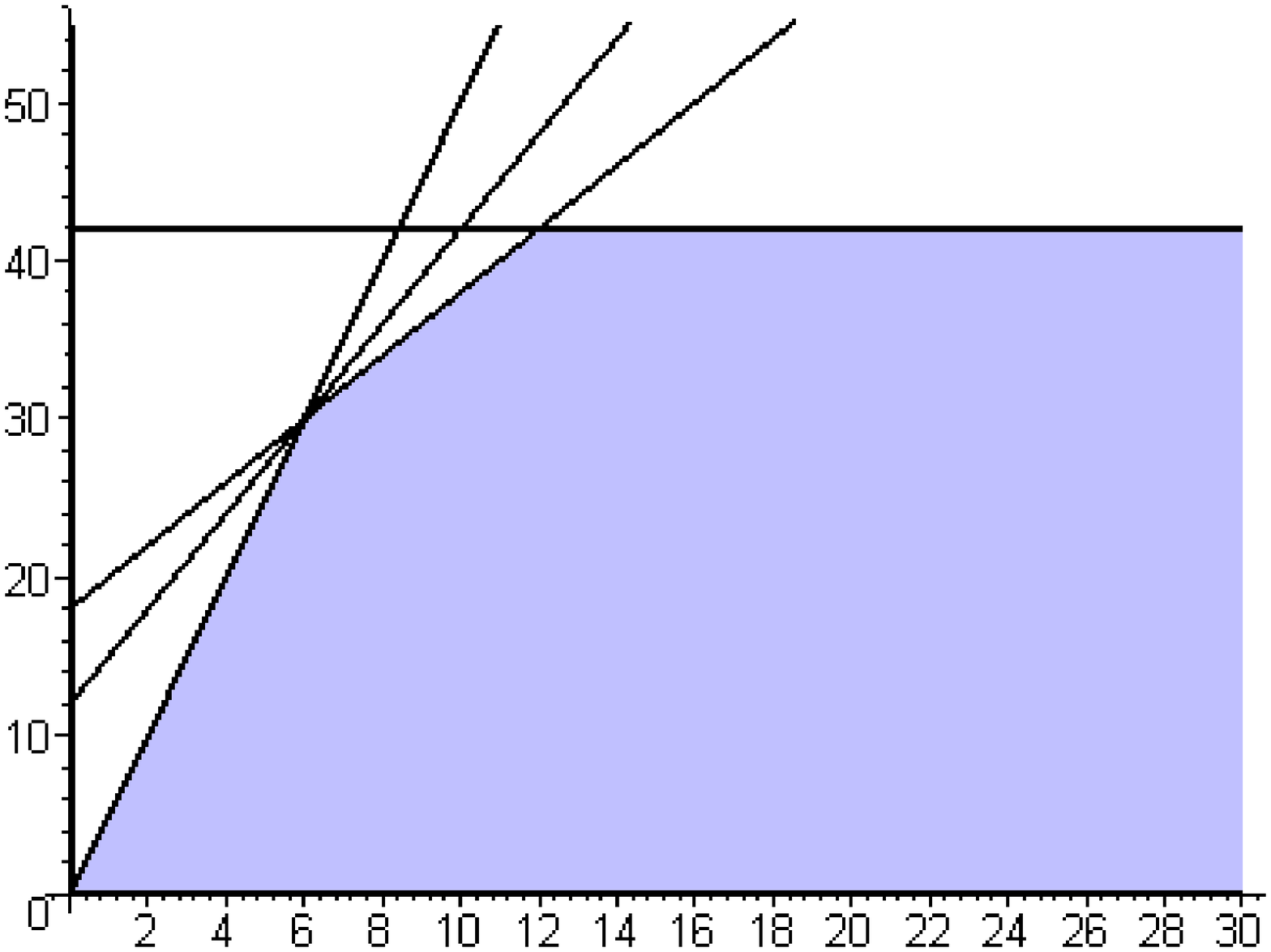}
\caption{$\Gamma ^{A}$ $\left( g\right) .$}
\end{figure}

The face functions are 
\begin{eqnarray*}
g_{(0,0)}(x,y) &=&(y^{3}-x^{2})^{5}, \\
\text{ }g_{(6,30)}(x,y)
&=&(y^{3}-x^{2})^{5}+(y^{3}-x^{2})^{3}x^{6}y^{3}+(y^{3}-x^{2})^{2}x^{12}, \\
g_{(12,40)}(x,y) &=&(y^{3}-x^{2})^{2}x^{12}+x^{24}.
\end{eqnarray*}
Since $g_{(6,30)}(x,y)=\left( y^{3}-x^{2}\right) ^{2}\left(
(y^{3}-x^{2})^{3}+(y^{3}-x^{2})+x^{12}\right) $ has singular points on $%
\left( K^{\times }\right) ^{2}$, $g(x,y)$\ is arithmetically \ degenerate.
\end{example}

\begin{example}
We set 
\begin{equation*}
h(x,y)=(y^{5}-x^{3})^{4}(y^{5}-ax^{3})(y^{5}-bx^{3})+x^{20}\in K[x,y],
\end{equation*}
with $1\neq a\neq b$.The polynomial \ $h(x,y)$ is a\ degenerate, in the
sense of Kouchnirenko,\ semi-quasihomogeneous polynomial with respect the
weight $(5,3)$. The origin of $K^{2}$ is a singular point of $h(x,y),$ and
there are no singularities on the torus $\left( K^{\times }\right) ^{2}$. In
this case the arithmetic Newton polygon of $h(x,y)$ is \ equal to $\Gamma
^{A}\left( h\right) =\{\Gamma _{h,1}$, $\Gamma _{h,a}$, $\Gamma _{h,b}\}$. \
The boundary of Newton polygon $\Gamma _{h,1}$ is conformed by the straight
segments 
\begin{eqnarray*}
w_{1,1}(z) &=&4z,\text{ \ }0\leq z\leq \frac{5}{2}, \\
w_{2,1}(z) &=&10,\text{ \ }z\geqslant \frac{5}{2},
\end{eqnarray*}
and the half-line \ $\left\{ \left( z,w\right) \in \mathbb{R}_{+}^{2}\mid
w=0\right\} $. The arithmetic \ Newton polygon $\Gamma _{h,1}$ is showed in
the next figure.
\bigskip

\begin{figure}[htp]
\hskip 4cm \epsfxsize=6cm \epsfbox{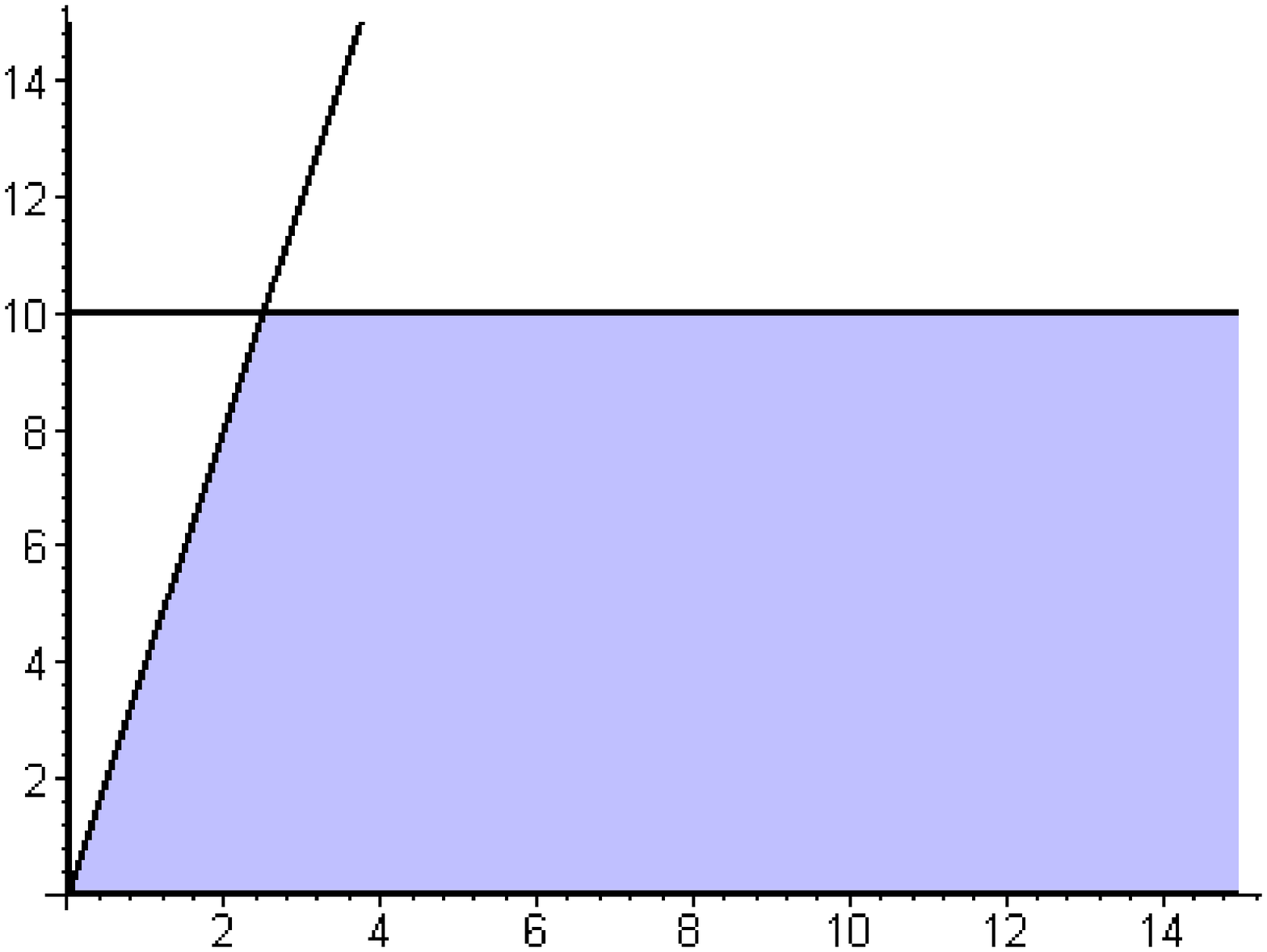}
\caption{$\Gamma _{h,1}.$}
\end{figure}

The face functions are 
\begin{eqnarray*}
h_{(0,0)}(x,y) &=&(y^{5}-x^{3})^{4}(y^{5}-ax^{3})(y^{5}-bx^{3}) \\
\text{ }h_{(\frac{5}{2},10)}(x,y)
&=&(y^{5}-x^{3})^{4}(y^{5}-ax^{3})(y^{5}-bx^{3})+x^{20}.
\end{eqnarray*}
Since 
\begin{equation*}
h_{(\frac{5}{2}%
,10)}(x,y)=(y^{5}-x^{3})^{4}(y^{5}-ax^{3})(y^{5}-bx^{3})+x^{20}
\end{equation*}
has no singular points on $\left( K^{\times }\right) $, $\ f(x,y)$\ is
arithmetically non \ degenerate with respect to $\Gamma _{h,1}$.

The boundary of the Newton polygon $\Gamma _{h,a}$ is conformed by the
straight segments 
\begin{eqnarray*}
w_{1,a}(z) &=&z,\text{ \ }0\leq z\leq 10, \\
w_{2,a}(z) &=&10,\text{ \ }z\geqslant 10,
\end{eqnarray*}
and the half-line \ $\left\{ \left( z,w\right) \in \mathbb{R}_{+}^{2}\mid
w=0\right\} $. The arithmetic \ Newton diagram $\Gamma _{h,a}$ is showed in
the next figure.

\begin{figure}[htp]
\hskip 4cm \epsfxsize=6cm \epsfbox{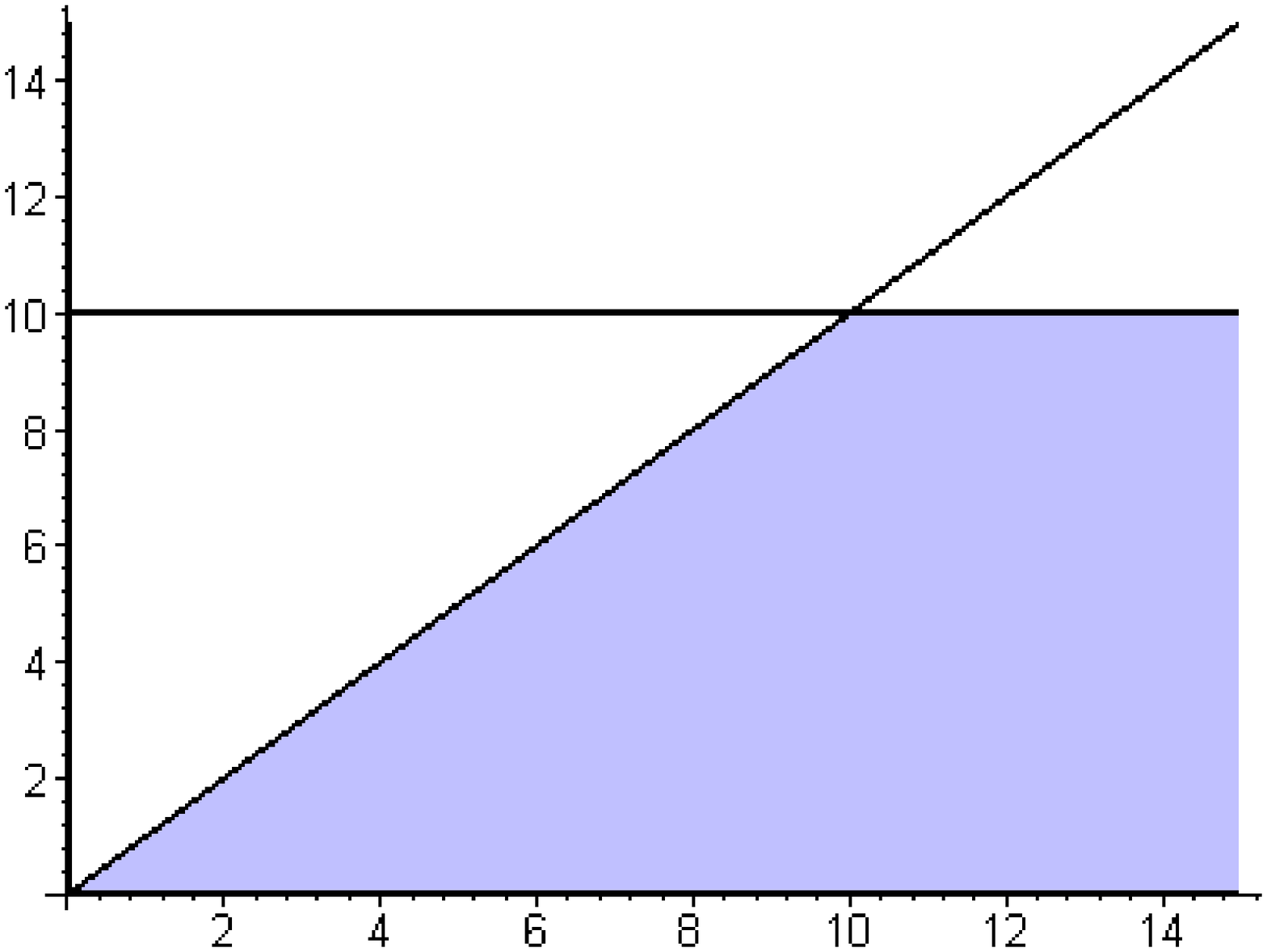}
\caption{$\Gamma _{h,a}.$}
\end{figure}

The face functions are 
\begin{eqnarray*}
h_{(0,0)}(x,y) &=&(y^{5}-x^{3})^{4}(y^{5}-ax^{3})(y^{5}-bx^{3}), \\
\text{ }h_{(10,10)}(x,y)
&=&(y^{5}-x^{3})^{4}(y^{5}-ax^{3})(y^{5}-bx^{3})+x^{20}.
\end{eqnarray*}
Since 
\begin{equation*}
h_{(10,10)}(x,y)=(y^{5}-x^{3})^{4}(y^{5}-ax^{3})(y^{5}-bx^{3})+x^{20}
\end{equation*}
has no singular points on $\left( K^{\times }\right) ^{2}$, $f(x,y)$\ is
arithmetically non \ degenerate with respect to $\Gamma _{h,a}$. The Newton
polygon $\Gamma _{h,b}$ is equal to $\Gamma _{h,a}$. Thus the polynomial $%
h(x,y)$ is arithmetically non degenerate with respect \ $\Gamma _{h}^{A}$.
\end{example}

\section{Degenerate integrals of type $Z(s,f,v,\Delta )$}

In this section using the notion of arithmetic non degeneracy introduced in
the above section, we \ shall \ compute explicitly degenerate \ integrals of
type 
\begin{equation*}
Z(s,f,v,\Delta )=\int\limits_{^{E\text{ }}}\mid f(x,y)\mid _{v}^{s}\mid
dxdy\mid ,
\end{equation*}
where 
\begin{equation*}
\Delta :=(a,b)\mathbb{R}_{+}\text{, and \ }E:=\{(x,y)\in O_{v}^{2}\mid
(v(x),v(y))\in \Delta \},
\end{equation*}
and $f(x,y)\in O_{v}[x,y]$ is\ a degenerate, in the sense of Kouchnirenko,
semi-quasi\-ho\-mo\-geneous polynomial with respect to the weight $(a,b)$,
but \ arithmetically \ non degenerate with respect $\Gamma ^{A}\left(
f\right) =\cup _{\theta }\Gamma _{f,\theta }$. From now on, we shall
consider only those $\Gamma _{f,\theta }$ for which $\theta \in O_{v}$. The
reason is that the $\Gamma _{f,\theta }$, with $\theta \in L_{v}\setminus
O_{v}$, do not contribute the poles of $Z(s,f,v,\Delta )$.

\begin{definition}
For a semi-quasihomogeneous polynomial $f(x,y)\in L_{v}[x,y]$ non degenerate
with respect to $\Gamma ^{A}(f)=\cup _{\left\{ \theta \ \in O_{v}\mid
f_{0}(1,\theta ^{a})=0\right\} }\Gamma _{f,\theta }$, we define 
\begin{equation}
\mathcal{P}(\Gamma _{f,\theta }):=\bigcup\limits_{i=1}^{r_{\theta }}\{-\frac{%
1}{\mathcal{E}}_{i},-\frac{(a+b)+\tau _{i}}{\mathcal{D}_{i+1}+\mathcal{E}%
_{i+1}\tau _{i}},-\frac{(a+b)+\tau _{i}}{\mathcal{D}_{i}+\mathcal{E}_{i}\tau
_{i}}\}\bigcup \bigcup_{\{\mathcal{E}_{r+1}\neq o\}}\{-\frac{1}{\mathcal{E}}%
_{r+1}\},
\end{equation}
and 
\begin{equation}
\mathcal{P}(\Gamma ^{A}(f)):=\bigcup\limits_{\{\theta \ \in O_{v}\mid
f_{0}(1,\theta ^{a})=0\}}\mathcal{P}(\Gamma _{f,\theta }).
\end{equation}
\end{definition}

The data \ $\mathcal{D}_{i}$, $\mathcal{E}_{i}$, $\tau _{i}$, for each $i$,
are obtained \ from the equations of the straight segments conforming the
boundary of $\Gamma _{f,\theta }$. \bigskip

The main result of this section is the following.

\begin{theorem}
\label{theorem1}Let $f(x,y)=\sum_{j=0}^{l_{f}}f_{j}(x,y)\in O_{v}[x,y]$ be a
semi-quasihomogeneous polynomial, with respect to the weight $(a,b)$, with $a
$, $b$ coprime, $f_{j}(x,y)$ as in (\ref{kqs}). If $f(x,y)$ is
arithmetically non degenerate with respect to $\Gamma ^{A}(f)$, then the
real parts of the poles of $Z(s,f,v,\Delta )$ belong to the set 
\begin{equation*}
\{-1\}\bigcup \{-\frac{a+b}{d_{0}}\}\bigcup \mathcal{P}(\Gamma ^{A}(f)).
\end{equation*}
\end{theorem}

The proof will be given at the end of this section, and it will be
accomplished \ by computing explicitly several $p-$adic integrals. Moreover
the proof provides an effective method \ for computing these integrals.

The following two examples illustrate the theorem.

\begin{example}
We set $f(x,y)=(y^{3}-x^{2})^{2}+x^{4}y^{4}\in L_{v}[x,y]$. We suppose that
the characteristic of the residue field is different from $2$. The
polynomial \ $f(x,y)$ is arithmetically \ non degenerate,
semi-quasihomogeneous polynomial with respect to the weight $(3,2)$ (see
example (\ref{example1}). In this case, we have the following data attached
\ the \ arithmetic Newton polygon\ $\Gamma ^{A}(f)=\Gamma _{f,1}$, $a=3$, $%
b=2$,$\mathcal{D}_{1}=d_{0}=12$, $\tau _{1}=4,\mathcal{E}_{1}=2$, $\mathcal{D%
}_{2}=20$. Thus according the above theorem, the real parts of the poles of
the integral $Z(s,f,v)$\ belong to the set 
\begin{equation*}
\{-1\}\bigcup \{-\frac{5}{12}\}\bigcup \{-\frac{1}{2},-\frac{9}{20}\}.
\end{equation*}
\end{example}

In subsection (\ref{sub1}), \ we computed explicitly the local zeta function 
$Z(s,f,v)$ (cf. (\ref{poles1})).

\begin{example}
We set $g(x,y)=(y^{3}-x^{2})^{2}(y^{3}-ax^{2})+x^{4}y^{4}\in L_{v}[x,y]$,
with $a%
%TCIMACRO{
%\TeXButton{notequiv}{\not\equiv%
%}}%
%BeginExpansion
\not\equiv%
%
%EndExpansion
1$ mod $\pi $. We assume that the characteristic of the residue field of $%
L_{v}$ is different from $2$. The polynomial $g(x,y)$ is quasihomogeneous
with respect to the weight $(3,2)$. The origin of $L_{v}^{2}$ is the only
singular \ point of \ $g(x,y)$.

The arithmetic Newton polygon $\Gamma ^{A}(g)=\{\Gamma _{g,1},$ $\Gamma
_{g,a}\}$. The boundary of the Newton polygon $\Gamma _{g,1}$ is conformed
by 
\begin{eqnarray*}
w_{1,1}(z)=2z,\text{ \ }0\leq z\leq 1, \\
w_{2,1}(z)=2,\text{ \ }z\geqslant 1, \\
w(z)=0,\text{ \ }z\geqslant 0.
\end{eqnarray*}
Thus $\mathcal{D}_{1}=d_{0}=18$, $\mathcal{E}_{1}=2$, $\mathcal{D}_{2}=20$, $%
\mathcal{E}_{2}=0$, $\tau _{0}=0$,$\ \tau _{1}=1$.The polynomial is
arithmetically non degenerate with respect $\Gamma _{g,1}$. According to the
previous theorem the contribution of the \ arithmetic Newton polygon $\Gamma
_{g,1}$ to the set of real parts of $Z(s,g,v)$\ is \ $\{-\frac{1}{2},-\frac{3%
}{10}\}.$

The boundary of the Newton polygon $\Gamma _{g,a}$ is conformed by 
\begin{eqnarray*}
w_{1,1}(z)=z,\text{ \ }0\leq z\leq 2, \\
w_{2,1}(z)=2,\text{ \ }z\geqslant 2. \\
w(z)=0,\text{ \ }z\geqslant 0.
\end{eqnarray*}
Thus $\mathcal{D}_{1}=d_{0}=18$, $\mathcal{E}_{1}=1$, $\mathcal{D}_{2}=20$, $%
\mathcal{E}_{2}=0$, $\tau _{0}=0$,$\ \tau _{1}=2$. The polynomial is
arithmetically non degenerate with respect $\Gamma _{g,a}.$ According to the
previous theorem the contribution of the \ arithmetic Newton polygon $\Gamma
_{g,a}$ to the set of real parts of $Z(s,g,v)$\ is \ $\{-1,$ $-\frac{7}{20}%
\}.$ Then\ according to the previous theorem the \ real parts of the poles
of $Z(s,g,v)$ belong to the set \ $\{-1\}\bigcup \{-\frac{5}{18}\}\bigcup \{-%
\frac{1}{2},-\frac{9}{20}\}\bigcup \{-\frac{7}{20}\}$. In subsection (\ref
{sub2}), we computed explicitly the local zeta $Z(s,g,v)$ (cf. (\ref{poles2}%
)).
\end{example}

\subsection{Some $p-$adic integrals}

\begin{proposition}
\label{prop1}Let $f(x,y)\in O_{v}[x,y]$ be a semi-quasihomogeneous
polynomial, with respect to the weight $(a,b)$, with $a$, $b$ coprime, and 
\begin{equation}
f^{(m)}(x,y):=\pi ^{-d_{0}m}f(\pi ^{am}x,\pi ^{bm}y)=\sum_{j=0}^{l_{f}}\pi
^{(d_{j}-d_{0})m}f_{j}(x,y),  \label{fm}
\end{equation}
with $f_{j}(x,y)$ as in (\ref{kqs}), and $m\geqq 1$. Then there exists a
measure-preserving \ bijection 
\begin{eqnarray*}
\Phi :O_{v}^{\times \text{ }2}\rightarrow O_{v}^{\times \text{ }2} \\
(u,w)\rightarrow (\Phi _{1}(x,y),\Phi _{2}(x,y))\text{,}
\end{eqnarray*}
such that $f^{(m)}\circ \Phi (x,y)=u^{N_{i}}w^{M_{i}}\widetilde{f^{(m)}}%
(u,w) $, with $\widetilde{f^{(m)}}(u,w)$ non vanishing identically \ on the
hypersurface $uw=0$, and 
\begin{equation}
\widetilde{f^{(m)}}(u,w)=\sum_{j=0}^{l_{f}}\pi ^{(d_{j}-d_{0})m}\widetilde{%
f_{j}}(u,w),  \label{fm1}
\end{equation}
with 
\begin{eqnarray}
\widetilde{f}_{j}(u,w)=c_{j}u^{A_{j}}w^{B_{j}}\prod_{i=1}^{l_{j}}\left(
w-\alpha _{i,j}\right) ^{e_{i,j}}\text{, \ \ or }  \notag \\
\widetilde{f}_{j}(u,w)=c_{j}u^{L_{j}}w^{F_{j}}\prod_{i=1}^{l_{j}}\left(
1-\alpha _{i,j}u\right) ^{e_{i,j}}.  \label{fm2}
\end{eqnarray}
\end{proposition}

\begin{proof}
We denote by $\Phi _{1}$ the map 
\begin{eqnarray*}
\ \Phi _{1}(u,w):O_{v}^{\times \text{ }2}\rightarrow O_{v}^{\times \text{ }2}
\\
(u,w)\rightarrow (x,y),
\end{eqnarray*}
\ with $x=u$, $y=uw$, and by $\Phi _{2}$ the map 
\begin{eqnarray*}
\ \Phi _{2}(u,w):O_{v}^{\times \text{ }2}\rightarrow O_{v}^{\times \text{ }2}
\\
(u,w)\rightarrow (x,y),
\end{eqnarray*}
with $x=uw$, $y=w$. \ 

Since \ $\mid \det \Phi _{1}^{\prime }(u,w)\mid _{v}=1$, and$\mid \det \Phi
_{2}^{\prime }(u,w)\mid _{v}=1$, the maps $\Phi _{1}$, $\Phi _{2}$\ \ give a
measure-preserving bijection of $O_{v}^{\times \text{ }2}$\ to itself. We
shall show that $\Phi $\ is a finite composition \ of maps of the form $\Phi
_{1}$ or $\Phi _{2}$. The proof will be accomplished by induction on \ $\min
\{a,b\}.$

\textbf{Case} \ $\min \{a,b\}=1$

In this case, it is sufficient to take $\Phi $ as follows:

\begin{equation*}
\Phi =\left\{ 
\begin{array}{ccc}
\underbrace{\Phi _{1}\circ \cdots \circ \Phi _{1}}_{b-\text{times}} & \text{%
if} & a=1\text{,} \\ 
\underbrace{\Phi _{2}\circ \cdots \circ \Phi _{2}}_{a-\text{times}} & \text{%
if} & b=1.
\end{array}
\right.
\end{equation*}

\textbf{Induction hypothesis}

Suppose by \ induction hypothesis that the proposition is valid for all
polynomials $f^{(m)}(x,y)$ of the form (\ref{fm}) satisfying \ 1$\leqq \min
\{a,b\}\leqq k$, with $k\geqq 1$.

\textbf{Case} \ $\min \{a,b\}=k+1,$ $k\geqq 2$

\qquad Let $f^{(m)}(x,y)$\ be a polynomial of the form (\ref{fm}) satisfying
\ $\min \{a,b\}=k+1$, $k\geqq 1$, and $a>b$. By applying the Euclidean
algorithm to $a$, $b$, we have that 
\begin{equation*}
a=q_{1}b+r_{1},\text{ \ \ \ \ \ \ \ \ }0\leqq r_{1}<b,
\end{equation*}

for some $q_{1}$, $r_{1}\in \mathbb{N}$. Because \ $a$, $b$ are coprime,
necessarily \ \ $1\leqq r_{1}<b.$

We set $\Psi = 
\begin{array}{c}
\underbrace{\Phi _{2}\circ \cdots \circ \Phi _{2}} \\ 
q_{1}-\text{times}
\end{array}
$, i.e. $x=uw^{q_{1}}$, $y=w,$ thus 
\begin{equation*}
f^{(m)}\circ \Psi (x,y)=u^{A_{i}}w^{B_{i}}f^{\ast (m)}(u,w)
\end{equation*}
with $\ f^{\ast (m)}(u,w)$ non vanishing identically on $uw=0$, and 
\begin{equation}
f^{\ast (m)}(u,w)=\sum_{j=0}^{l_{f}}\pi ^{(d_{j}-d_{0})m}f_{j}^{\ast }(u,w)%
\text{, }  \label{fm3}
\end{equation}
with 
\begin{equation*}
f_{j}^{\ast }(u,w)=c_{j}u^{C_{j}}w^{D_{j}}\prod_{i=1}^{l_{j}}\left(
w^{r_{1}}-\alpha _{i,j}u^{b}\right) ^{e_{i,j}}\text{.}
\end{equation*}

Since $\min \{r_{1},b\}=r_{1},$ and $1\leqq r_{1}\leqq b-1=k$, by applying
the induction hypothesis applied to \ the polynomial \ $f^{\ast (m)}(u,v)$
in (\ref{fm3}), there exists a map $\Theta $, that is a finite composition
of maps of the form $\Phi _{1}$ or $\Phi _{2}$, such that $\left(
f^{(m)}\circ \Psi \right) \circ \Theta =f^{(m)}\circ \left( \Psi \circ
\Theta \right) =u^{N_{i}}w^{M_{i}}$\ $\widetilde{f^{(m)}}(u,v)$, and $%
\widetilde{f^{(m)}}(u,v)$ has \ all announced properties in the proposition.
Therefore, it is sufficient to take $\Phi =\Psi \circ \Theta $.

The case $b>a$ is proved in an analogous form.
\end{proof}

\begin{remark}
\label{remark1}We can suppose without loss of generality that 
\begin{equation}
\widetilde{f^{(m)}}(u,w)=\sum_{j=0}^{l_{f}}\pi ^{(d_{j}-d_{0})m}\widetilde{%
f_{j}}(u,w),
\end{equation}
with 
\begin{equation}
\widetilde{f}_{j}(u,w)=c_{j}u^{A_{j}}w^{B_{j}}\prod_{i=1}^{l_{j}}\left(
w-\alpha _{i,j}\right) ^{e_{i,j}},\text{ \ with }c_{j}\text{, }\alpha
_{i,j}\in L_{v}^{\times }\text{.}  \label{fm3a}
\end{equation}
\end{remark}

We define 
\begin{equation*}
I(s,f^{(m)},v,O_{v}^{\times 2}):=\int\limits_{^{O_{v}^{\times 2}\text{ }%
}}\mid f^{(m)}(x,y)\mid _{v}^{s}\mid dxdy\mid .
\end{equation*}

\begin{proposition}
\label{prop2}Let $f(x,y)\in O_{v}[x,y]$ be a semi-quasihomogeneous
polynomial, with respect to the weight $(a,b)$, with $a$, $b$ coprime, and 
\begin{equation}
f^{(m)}(x,y)=\pi ^{-d_{0}m}f(\pi ^{am}x,\pi ^{bm}y)=\sum_{j=0}^{l_{f}}\pi
^{(d_{j}-d_{0})m}f_{j}(x,y),\text{ \ \ }m\geqq 1,
\end{equation}
with $f_{j}(x,y)$ as in (\ref{kqs}). There exists a constant $M_{0}$ such
that if $m\geqslant M_{0}$, then 
\begin{eqnarray}
I(s,f^{(m)},v,O_{v}^{\times 2}) &=&U_{0}(q^{_{-s}})  \notag \\
+\sum\limits_{\{\theta \in O_{v}\mid f_{0}(1,\theta
^{a})=0\}}\sum\limits_{k=1}^{\infty }q^{-k}\int\limits_{O_{v}^{\times \text{ 
}}\times O_{v}^{\times \text{ }}\text{ }} &\mid &\widetilde{f^{(m)}}%
(x,\theta +\pi ^{k}y)\mid _{v}^{s}\mid dxdy\mid ,
\end{eqnarray}
where $U_{0}(q^{_{-s}})$\ is a polynomial with rational coefficients.
Moreover, the polynomial $U_{0}(q^{-s})$ and the constant $M_{0}$ can be
computed effectively.
\end{proposition}

\begin{proof}
Proposition \ref{prop1} implies that

\begin{equation}
I(s,f^{(m)},v,O_{v}^{\times 2})=\int\limits_{O_{v}^{\times 2}}\mid
f^{(m)}(x,y)\mid _{v}^{s}\mid dxdy\mid =\int\limits_{O_{v}^{\times 2}}\mid 
\widetilde{f^{(m)}}(x,y)\mid _{v}^{s}\mid dxdy\mid .  \label{I1}
\end{equation}

We set 
\begin{equation*}
R\left( f_{0}\right) :=\left\{ \theta \in O_{v}\mid f_{0}\left( 1,\theta
^{a}\right) =0\right\} ,
\end{equation*}
\begin{equation*}
l\left( f_{0}\right) :=\max_{
\begin{array}{c}
\theta \neq \theta ^{\prime } \\ 
\theta ,\theta ^{\prime }\in R\left( f_{0}\right) 
\end{array}
}\left\{ v\left( \theta -\theta ^{\prime }\right) \right\} ,
\end{equation*}
and 
\begin{equation*}
B(\theta )=B\left( l\left( f_{0}\right) ,\theta \right) :=O_{v}^{\times
}\times \left( \theta +\pi ^{1+l\left( f_{0}\right) }O_{v}\right) \text{, \
for }\theta \in O_{v},\text{with }v\left( \theta \right) \leq l\left(
f_{0}\right) .
\end{equation*}

With the above notation, if $\theta =\alpha _{i_{0},0}$, then 
\begin{equation*}
v\left( \widetilde{f}_{j}\left( x,y\right) \right) =\left\{ 
\begin{array}{cc}
v(c_{j})+B_{j}v(\theta )+\sum\limits_{i=1}^{l_{j}}e_{i,j}v\left( \theta
-\alpha _{i,j}\right)  & \text{if \ }f_{j}\left( 1,\theta \right) \neq 0, \\ 
\begin{array}{c}
v(c_{j})+B_{j}v(\theta )+\sum\limits_{
\begin{array}{c}
i=1 \\ 
i\neq i_{0}
\end{array}
}^{l_{j}}e_{i,j}v\left( \theta -\alpha _{i,j}\right)  \\ 
+e_{i_{0},j}\left( 1+l\left( f_{0}\right) \right) +e_{i_{0},j}v\left( \frac{%
y-\theta }{\pi ^{1+l\left( f_{0}\right) }}\right) 
\end{array}
& \text{if \ }f_{j}\left( 1,\theta \right) =0,
\end{array}
\right. 
\end{equation*}

for every $\left( x,y\right) \in B\left( l\left( f_{0}\right) ,\theta
\right) $. We \ put 
\begin{equation*}
Const\left( j,\theta \right) :=\left\{ 
\begin{array}{cc}
v(c_{j})+B_{j}v(\theta )+\sum\limits_{i=1}^{l_{j}}e_{i,j}v\left( \theta
-\alpha _{i,j}\right)  & \text{if \ }f_{j}\left( 1,\theta \right) \neq 0, \\ 
\begin{array}{c}
v(c_{j})+B_{j}v(\theta )+\sum\limits_{
\begin{array}{c}
i=1 \\ 
i\neq i_{0}
\end{array}
}^{l_{j}}e_{i,j}v\left( \theta -\alpha _{i,j}\right)  \\ 
+e_{i_{0},j}\left( 1+l\left( f_{0}\right) \right) 
\end{array}
& \text{if \ }f_{j}\left( 1,\theta \right) =0.
\end{array}
\right. 
\end{equation*}

Then, if $\theta \notin R\left( f_{0}\right) $, and $m\geqq M_{0}$, with 
\begin{equation*}
M_{0}:=\max_{
\begin{array}{c}
1\leq j\leq l_{f} \\ 
\theta 
\end{array}
}\left\{ Const\left( 0,\theta \right) +1-Const\left( j,\theta \right)
\right\} ,
\end{equation*}
it \ holds that 
\begin{equation*}
v\left( \widetilde{f^{(m)}}\left( x,y\right) \right) =Const\left( 0,\theta
\right) \text{, \ }
\end{equation*}
for every $\left( x,y\right) \in B\left( l\left( f_{0}\right) ,\theta
\right) $. From the above calculations follow that 
\begin{equation}
\int\limits_{B\left( l\left( f_{0}\right) ,\theta \right) }\mid \widetilde{%
f^{(m)}}(x,y)\mid _{v}^{s}\mid dxdy\mid =\mid Const\left( 0,\theta \right)
\mid _{v}^{s}\text{, if \ }\theta \notin R\left( f_{0}\right) ,\text{ and }%
m\geqq M_{0}.  \label{I1.1}
\end{equation}
By subdividing $O_{v}^{\times }\times O_{v}^{\times }$\ into equivalence
classes modulo $\pi ^{1+l\left( f_{0}\right) }$, and using (\ref{I1}), and (%
\ref{I1.1}), \ we have \ that

\begin{equation}
I(s,f^{(m)},v,O_{v}^{\times 2})=U_{0}(q^{-s})+\sum\limits_{\theta \in
R\left( f_{0}\right) }\int\limits_{B\left( \theta \right) \text{ }}\mid 
\widetilde{f^{(m)}}(x,y)\mid _{v}^{s}\mid dxdy\mid ,  \label{I2}
\end{equation}
for $m\geqq M_{0}$, where $U_{0}(q^{-s})$ is a polynomial with rational
coefficients. Moreover, the polynomial $U_{0}(q^{-s})$ and the constant $%
M_{0}$ can be computed effectively.
\end{proof}

The integral $\int_{B\left( \theta \right) }\mid \widetilde{f^{(m)}}%
(x,y)\mid _{v}^{s}\mid dxdy\mid ,$ for $\theta \in R\left( f_{0}\right) $, \
admits the following expansion 
\begin{eqnarray}
\sum\limits_{k=1}^{\infty }q^{-k}\int\limits_{O_{v}^{\times \text{ }}\times
O_{v}^{\times \text{ }}\text{ }} &\mid &\widetilde{f^{(m)}}(x,\theta +\pi
^{k}y)\mid _{v}^{s}\mid dxdy\mid =  \notag \\
\sum\limits_{k=1}^{^{1+l\left( f_{0}\right)
}}q^{-k}\int\limits_{O_{v}^{\times \text{ }}\times O_{v}^{\times \text{ }}%
\text{ }} &\mid &\widetilde{f^{(m)}}(x,\theta +\pi ^{k}y)\mid _{v}^{s}\mid
dxdy\mid  \notag \\
+\sum\limits_{k=2+l(f_{0})}^{\infty }q^{-k}\int\limits_{O_{v}^{\times \text{ 
}}\times O_{v}^{\times \text{ }}\text{ }} &\mid &\widetilde{f^{(m)}}%
(x,\theta +\pi ^{k}y)\mid _{v}^{s}\mid dxdy\mid =  \notag \\
\frac{U_{1}(q^{-s})}{1-q^{-1-s}}+\sum\limits_{k=2+l(f_{0})}^{\infty
}q^{-k}\int\limits_{O_{v}^{\times \text{ }}\times O_{v}^{\times \text{ }}%
\text{ }} &\mid &\widetilde{f^{(m)}}(x,\theta +\pi ^{k}y)\mid _{v}^{s}\mid
dxdy\mid ,  \label{fm0}
\end{eqnarray}
(cf. lemma \ref{lemma2}). We set

\begin{equation*}
J(s,m,\theta ):=\sum\limits_{k=2+l(f_{0})}^{\infty
}q^{-k}\int\limits_{O_{v}^{\times \text{ }}\times O_{v}^{\times \text{ }}%
\text{ }}\mid \widetilde{f^{(m)}}(x,\theta +\pi ^{k}y)\mid _{v}^{s}\mid
dxdy\mid .
\end{equation*}
Then from (\ref{I2}) and (\ref{fm0}) follow that

\begin{equation}
I(s,f^{(m)},v,O_{v}^{\times 2})=\frac{U(q^{-s})}{1-q^{-1-s}}%
+\sum\limits_{\theta \in R\left( f_{0}\right) }J(s,m,\theta ),\text{ \ }%
m\geqslant M_{0},  \label{fm0a}
\end{equation}
where \ $U(q^{-s})$ is a polynomial with rational coefficients.

The polynomial $\widetilde{f^{(m)}}(x,\theta +\pi ^{k}y)$\ can expressed as 
\begin{equation}
\widetilde{f^{(m)}}(x,\theta +\pi ^{k}y)=\sum_{j=0}^{l_{f}}C(j,\theta )\pi
^{(d_{j}-d_{0})m+ke_{j,\theta }}\gamma _{j}(x,y)y^{e_{j,\theta }},
\label{fml}
\end{equation}
for\ \ $k\geqq 2+l\left( f_{0}\right) $, with $C(j,\theta )=Const(j,\theta
)c_{j}\in L_{v}^{\times }$, and $\left| \gamma _{j}(x,y)\right| _{v}=1$, for
every $(x,y)\in O_{v}^{\times }\times O_{v}^{\times }$.\ 

Given a real number $x$, \ $[x]$ denotes the greatest integer less than or
equal to $x$.

\begin{proposition}
\label{prop3}With the hypothesis of proposition \ref{prop2}. If $f(x,y)$ is
arithmetically non degenerate with respect to $\Gamma _{f,\theta \text{,}}$
then

\begin{eqnarray}
J(s,m,\theta ) &=&(1-q^{-1})^{2}\sum_{i=0}^{r-1}\frac{q^{-(1+s\mathcal{E}%
_{i+1})}}{1-q^{-(1+s\mathcal{E}_{i+1})}}q^{-(\mathcal{D}_{i+1}-d_{0})ms-(1+s%
\mathcal{E}_{i+1})[m\tau _{i}]}  \notag \\
&&-(1-q^{-1})^{2}\sum_{i=0}^{r-1}\frac{1}{1-q^{-(1+s\mathcal{E}_{i+1})}}q^{-(%
\mathcal{D}_{i+1}-d_{0})ms-(1+s\mathcal{E}_{i+1})[m\tau _{i+1}]}  \notag \\
&&+(1-q^{-1})^{2}\frac{q^{-(1+s\mathcal{E}_{r+1})}}{1-q^{-(1+s\mathcal{E}%
_{r+1})}}q^{-(\mathcal{D}_{r+1}-d_{0})ms-(1+s\mathcal{E}_{r+1})[m\tau _{r}]}
\notag \\
&&+\sum_{i=1}^{r}q^{-(\mathcal{D}_{i}-d_{0})ms-(1+s\mathcal{E}_{i})[m\tau
_{i}]}I_{i}(s),  \label{I(s,fm)}
\end{eqnarray}
where $\tau _{i}$, $i=0,1,\cdots ,r$, are the abscissas of the vertices of \ 
$\Gamma _{f,\alpha _{j,0\text{ }}}$, and $I_{i}(s)=\frac{M_{i}(q^{-s},m)}{%
1-q^{-1-s}}$, $M_{i}(q^{-s},m)\in \mathbb{Q}[q^{-s}],$ $i=1,2,\cdots ,r$.
Moreover, there exists a \ constant $c(f)$ such that if $m\geqq c(f)$, then $%
I_{i}(s)$ does not dependent on $m$.
\end{proposition}

\begin{proof}
The proof \ is based on an explicit computation of the $v-$adic absolute
value of \ $\widetilde{f^{(m)}}(x,\theta +\pi ^{k}y)$, when $x$, $y\in
O_{v}^{\times }$. In order to accomplish it, we attach to $\widetilde{f^{(m)}%
}(x,\theta +\pi ^{k}y)$ the convex set $\Gamma _{\widetilde{f^{(m)}}%
(x,\theta +\pi ^{k}y)}$ defined as follows. We associate to each term 
\begin{equation}
C(j,\theta )\pi ^{(d_{j}-d_{0})m+ke_{j,\theta }}\gamma
_{j}(x,y)y^{e_{j,\theta }}  \label{cor1}
\end{equation}
of \ $\widetilde{f^{(m)}}(x,\theta +\pi ^{k}y)$ (see (\ref{fml})\ )\ an
straight line of the form 
\begin{equation*}
\widetilde{w}_{j,\theta }(\widetilde{z}):=(d_{j}-d_{0})m+e_{j,\theta }%
\widetilde{z},\text{ \ }j=0,1,\cdots ,l_{f},
\end{equation*}
and associate to$\widetilde{\text{ }f^{(m)}}(x,\theta +\pi ^{k}y)$ the
convex set 
\begin{equation}
\Gamma _{\widetilde{\text{ }f^{(m)}}(x,\theta +\pi ^{k}y)\mathit{\ }}=\{(%
\widetilde{z},\widetilde{w})\in \mathbb{R}_{+}^{2}\mid \widetilde{w}\leqq
\min_{0\leq j\leq l_{f\text{ }}}\{\widetilde{w}_{j,\theta }(\widetilde{z}%
)\}\}.
\end{equation}

Now, we set 
\begin{equation*}
\begin{array}{cccc}
\Omega _{m}: & \mathbb{R}^{2} & \longrightarrow  & \mathbb{R}^{2} \\ 
& (u,w) & \longrightarrow  & (\widetilde{u},\widetilde{w})\text{,}
\end{array}
\end{equation*}
with $u=\frac{\widetilde{u}}{m},$ $z=\frac{\widetilde{u}}{m}$. \ Then $%
\Omega _{m}(\Gamma _{f,\theta })=\Gamma _{\widetilde{f^{(m)}}(x,\theta +\pi
^{k}y)}$ (see definition \ref{4.1}). The homomorphism $\Omega _{m}$ sends
the topological boundary of $\Gamma _{f,\theta }$ into the topological
boundary of $\Gamma _{\widetilde{f^{(m)}}(x,\theta +\pi ^{j}y)}$. Thus the
vertices of $\Gamma _{\widetilde{f^{(m)}}(x,\theta +\pi ^{k}y)}$ are

\begin{equation}
\Omega _{m}(Q_{i}):=\left\{ 
\begin{array}{cc}
(0,0) & \text{if \ \ }i=0, \\ 
(m\tau _{i},(\mathcal{D}_{i}-d_{0})m+m\mathcal{E}_{i}\tau _{i}) & \text{if \ 
}i=1,2,\cdots ,r\text{,}
\end{array}
\right. 
\end{equation}
where the $\tau _{i}\ $\ are the abscissas of the \ vertices of $\Gamma
_{f^{(m)},\theta }$ (cf. (\ref{new3}), (\ref{new4})).

If 
\begin{equation*}
m\geqq \max_{j,\theta }\left\{ v\left( C(j,\theta )\right) \right\} ,
\end{equation*}
then $\mid \widetilde{f^{(m)}}(x,\theta +\pi ^{k}y)\mid _{v}$ \ can be
computed from $\Gamma _{\widetilde{f^{(m)}}(x,\theta +\pi ^{k}y)}$, as
follows (we use the notation \ introduced at (\ref{new2}), (\ref{new3}), and
(\ref{new4})):

If $m\tau _{i}<k<m\tau _{i+1}$, $i=0,1,\cdots ,r-1$, then 
\begin{equation}
\mid \widetilde{f^{(m)}}(x,\theta +\pi ^{k}y)\mid _{v}=q^{-(\mathcal{D}%
_{i+1}-d_{0})m-\mathcal{E}_{i+1}k\ }\left| C(i+1,\theta )\right| _{v}\text{, 
}(x,y)\in O_{v}^{\times \text{ }}\times O_{v}^{\times \text{ }}.
\label{case1}
\end{equation}

If $k>m\tau _{r}$, then 
\begin{equation}
\mid \widetilde{f^{(m)}}(x,\theta +\pi ^{k}y)\mid _{v}=q^{-(\mathcal{D}%
_{r+1}-d_{0})m-\mathcal{E}_{r+1}k\ }\left| C(r+1,\theta )\right| _{v}\text{, 
}(x,y)\in O_{v}^{\times \text{ }}\times O_{v}^{\times \text{ }}.
\label{case2}
\end{equation}

If \ $k=m\tau _{i}$, $i=1,2,\cdots ,r$, then 
\begin{eqnarray}
&\mid &\widetilde{f^{(m)}}(x,\theta +\pi ^{k}y)\mid _{v}=  \label{case3} \\
q^{-(\mathcal{D}_{i}-d_{0})m-\mathcal{E}_{i}k} &\mid &\widetilde{f_{\Omega
_{m}(Q_{i})}^{(m)}}(x,y)+\pi ^{m\left( \mathcal{D}_{i+1}-\mathcal{D}%
_{i}\right) }\left( \cdots \right) \mid _{v}\text{, }(x,y)\in O_{v}^{\times 
\text{ }}\times O_{v}^{\times \text{ }}\text{,}  \notag
\end{eqnarray}

with

\begin{equation*}
\widetilde{f_{\Omega _{m}(Q_{i})}^{(m)}}(x,y)+\pi ^{m\left( \mathcal{D}%
_{i+1}-\mathcal{D}_{i}\right) }\left( \cdots \right) =\widetilde{f_{\Omega
_{m}(Q_{i})}^{(m)}}(x,y)+
\end{equation*}
\ 
\begin{equation*}
\pi ^{m\left( \mathcal{D}_{i+1}-\mathcal{D}_{i}\right) }\left( \text{terms
with weighted degree }\geqq \mathcal{D}_{i+1}\right) ,
\end{equation*}
and 
\begin{equation*}
\ \widetilde{f_{_{\Omega _{m}(Q_{i})}}^{(m)}}(x,y)=\sum\limits_{\widetilde{w}%
_{i,\theta }(\Omega _{m}(Q_{i}))=0}\gamma _{i}(x,y)y^{e_{i,\theta }}\text{,}
\end{equation*}
here $\widetilde{w}_{i,\theta }(\widetilde{z})$ is the straight line
corresponding to the term 
\begin{equation*}
\pi ^{(d_{i}-d_{0})m+ke_{i,\theta }}\gamma _{i}(x,y)y^{e_{i,\theta }}.
\end{equation*}

In addition, we note any $k\in N$, $k\geqq 1$ satisfies only one \ of the
following conditions:

\begin{equation*}
\left[ m\tau _{i}\right] +1\leq k\leq \left[ m\tau _{i+1}\right] -1,\text{ \ 
}i=0,1,\cdots ,r-1;
\end{equation*}
\begin{equation*}
k=\left[ m\tau _{i}\right] ,\text{ \ }i=0,1,\cdots ,r;
\end{equation*}

or 
\begin{equation*}
k\geqq \left[ m\tau _{r}\right] +1.
\end{equation*}

Then by using (\ref{case1}), (\ref{case2}), (\ref{case3}), and the previous
observation, it follows that

\begin{eqnarray}
J(s,m,\theta ) &=&\sum\limits_{k=2+l(f_{0})}^{\infty
}q^{-k}\int\limits_{O_{v}^{\times \text{ }}\times O_{v}^{\times \text{ }}%
\text{ }}\mid \widetilde{f^{(m)}}(x,\theta +\pi ^{k}y)\mid _{v}^{s}\mid
dxdy\mid =  \notag \\
&&(1-q^{-1})^{2}\sum_{i=0}^{r-1}\left( q^{-(\mathcal{D}_{i+1}-d_{0})ms\
}\sum_{j=[m\tau _{i}]+1}^{[m\tau _{i+1}]-1}q^{-(1+s\mathcal{E}%
_{i+1})j}\right)   \notag \\
&&+(1-q^{-1})^{2}\sum_{j=[m\tau _{r}]+1}^{\infty }q^{-(\mathcal{D}%
_{r+1}-d_{0})ms-(1+s\mathcal{E}_{r+1})j}  \notag \\
&&+\sum_{i=1}^{r}q^{-(\mathcal{D}_{i}-d_{0})ms-(1+s\mathcal{E}_{i})[m\tau
_{i}]}I_{i}(s),  \label{I(s,fm)2}
\end{eqnarray}

with 
\begin{equation*}
I_{i}(s):=\int\limits_{O_{v}^{\times 2}}\mid \widetilde{f_{\Omega
_{m}(Q_{i})}^{(m)}}(x,y)+\pi ^{m\left( \mathcal{D}_{i+1}-\mathcal{D}%
_{i}\right) }\left( \text{higher order terms }\right) \mid _{v}\mid dxdy\mid
.
\end{equation*}
Since $\widetilde{f_{\Omega _{m}(Q_{i})}^{(m)}}(x,y)+\pi ^{\left( \mathcal{D}%
_{i+1}-\mathcal{D}_{i}\right) }\left( \text{higher order terms }\right) $,
and $\widetilde{f_{\Omega _{m}(Q_{i})}^{(m)}}(x,y)$ do not have singular
points on $\left( L_{v}^{\times }\right) ^{2}$(cf. proposition \ref{prop1}),
\ there exists a constant $c_{0}(f)$ such that 
\begin{equation}
I_{i}(s)=\int\limits_{O_{v}^{\times 2}}\mid \widetilde{f_{\Omega
_{m}(Q_{i})}^{(m)}}(x,y)\mid _{v}\mid dxdy\mid =\frac{U_{i}(q^{-s})}{%
1-q^{-1-s}}\text{, \ \ }U_{i}(q^{-s})\in \mathbb{Q}[q^{-s}]\text{, }%
i=1,\ldots ,r,  \label{Ii}
\end{equation}
for $m\geqq c_{0}(f)$ (cf. lemma \ref{lemma2}). Finally, the result follows
from (\ref{I(s,fm)2}), by using the algebraic identity 
\begin{equation}
\sum\limits_{k=A}^{B}z^{k}=\frac{z^{A}-z^{B+1}}{1-z}.  \label{identity}
\end{equation}

The constant $c(f)=\max \left\{ c_{0}(f),\max_{j,\theta }\left\{ v\left(
C(j,\theta )\right) \right\} \right\} .$
\end{proof}

\subsection{Some algebraic identities\label{saidentities}}

The following algebraic identities follows easily from (\ref{identity}).
These identities will be used later on for the explicit computation of
certain $p-$adic integrals. Let $\ s$ be a complex number, and $i$, $c_{0}$
non-negative integers , with $\ c_{0}\geqq 1$.

\begin{eqnarray*}
S_{1}(s,i,c_{0}) &:&=\sum_{m=c_{0}}^{\infty }q^{-(a+b)m-d_{0}ms-(\mathcal{D}%
_{i+1}-d_{0})ms-(1+s\mathcal{E}_{i+1})[m\tau _{i}]} \\
&=&\left\{ 
\begin{array}{ccc}
{\frac{q^{-c_{0}(a+b)-c_{0}d_{0}s}}{1-q^{-(a+b)-d_{0}s}}} & \text{if} & i=0,
\\ 
\sum_{l=0}^{\mathcal{E}_{i}-\mathcal{E}_{i+1}-1}\left( \frac{%
q^{-B_{i,l}-A_{l}(\gamma _{i}+\delta _{i}s)}}{1-q^{-\gamma _{i}-\delta _{i}s}%
}\right) & \text{if} & \text{\ }i=1,\ldots ,r,
\end{array}
\right.
\end{eqnarray*}

with 
\begin{equation*}
\begin{tabular}{|l|}
\hline
$\tau _{i}=\frac{\mathcal{D}_{_{i+1}}-\mathcal{D}_{i}}{\mathcal{E}_{i}-%
\mathcal{E}_{_{i+1}}}$,$\text{ }$ \\ \hline
$A_{l}:=\text{\ }\left\{ 
\begin{array}{ccc}
\left[ \frac{c_{0}}{\mathcal{E}_{i}-\mathcal{E}_{_{i+1}}}\right] & \text{ if 
} & l\neq 0\text{,} \\ 
1 & \text{ if } & l=0\text{, and }\left[ \frac{c_{0}}{\mathcal{E}_{i}-%
\mathcal{E}_{_{i+1}}}\right] =0\text{,} \\ 
\left[ \frac{c_{0}}{\mathcal{E}_{i}-\mathcal{E}_{_{i+1}}}\right] & \text{ if 
} & l=0\text{, and }\left[ \frac{c_{0}}{\mathcal{E}_{i}-\mathcal{E}_{_{i+1}}}%
\right] \geqq 1\text{,}
\end{array}
\right. $ \\ \hline
$B_{i,l}:=(a+b)l+[l\tau _{i}]+s([l\tau _{i}]\mathcal{E}_{_{i+1}}+l\mathcal{D}%
_{_{i+1}})$, \\ \hline
$\gamma _{i}:=(a+b)(\mathcal{E}_{i}-\mathcal{E}_{i+1})+(\mathcal{D}_{i+1}-%
\mathcal{D}_{i})$, \\ \hline
$\delta _{i}:=\mathcal{D}_{i+1}(\mathcal{E}_{i}-\mathcal{E}_{_{i+1}})+(%
\mathcal{D}_{_{i+1}}-\mathcal{D}_{i})\mathcal{E}_{_{i+1}}$. \\ \hline
\end{tabular}
\end{equation*}

\begin{eqnarray*}
S_{2}(s,i,c_{0}) &:&=\sum_{m=1}^{\infty }q^{-(a+b)m-d_{0}ms-(\mathcal{D}%
_{i+1}-d_{0})ms-(1+s\mathcal{E}_{i+1})[m\tau _{i+1}]} \\
&=&\sum_{l=0}^{\mathcal{E}_{i+1}-\mathcal{E}_{i+2}-1}\left( \frac{%
q^{-G_{i,l}-H_{l}(\rho _{i}+\sigma _{i}s)}}{1-q^{-\rho _{i}-\sigma _{i}s}}%
\right) ,\text{ \ }i=0,\cdots ,r-1,
\end{eqnarray*}
with 
\begin{equation*}
\begin{tabular}{|l|}
\hline
$\tau _{i+1}=\frac{\mathcal{D}_{_{i+2}}-\mathcal{D}_{i+1}}{\mathcal{E}_{i+1}-%
\mathcal{E}_{i+2}}$,$\text{ }$ \\ \hline
$H_{l}:=\text{\ }\left\{ 
\begin{array}{ccc}
\left[ \frac{c_{0}}{\mathcal{E}_{i+1}-\mathcal{E}_{i+2}}\right] & \text{ if }
& l\neq 0\text{,} \\ 
1 & \text{ if } & l=0\text{, and }\left[ \frac{c_{0}}{\mathcal{E}_{i+1}-%
\mathcal{E}_{i+2}}\right] =0\text{,} \\ 
\left[ \frac{c_{0}}{\mathcal{E}_{i+1}-\mathcal{E}_{i+2}}\right] & \text{ if }
& l=0\text{, and }\left[ \frac{c_{0}}{\mathcal{E}_{i+1}-\mathcal{E}_{i+2}}%
\right] \geqq 1\text{,}
\end{array}
\right. $ \\ \hline
$G_{i,l}:=(a+b)l+[l\tau _{i+1}]+s([l\tau _{i+1}]\mathcal{E}_{_{i+1}}+l%
\mathcal{D}_{_{i+1}})$, \\ \hline
$\rho _{i}:=(a+b)(\mathcal{E}_{i+1}-\mathcal{E}_{i+2})+(\mathcal{D}_{i+2}-%
\mathcal{D}_{i+1})$, \\ \hline
$\delta _{i}:=\mathcal{D}_{i+1}(\mathcal{E}_{i+1}-\mathcal{E}_{_{i+2}})+(%
\mathcal{D}_{_{i+2}}-\mathcal{D}_{i+1})\mathcal{E}_{_{i+1}}$. \\ \hline
\end{tabular}
\end{equation*}

\begin{remark}
\label{remmark2}With the above notation, the real parts of the poles of \ $%
S_{1}(s,i,c_{0})$, $S_{2}(s,i,c_{0})$, \ belong to the set \ $P(\Gamma
_{f,\theta })\cup \{-\frac{a+b}{d_{0}}\}.$\textrm{\ }
\end{remark}

\subsection{Proof of Theorem 5.1}

The integral $Z(s,f,v,\Delta )$ admits the following expansion:

\begin{eqnarray}
Z(s,f,v,\Delta ) &=&\sum\limits_{m=1}^{\infty }\int\limits_{\pi
^{am}O_{v}^{\times }\times \pi ^{bm}O_{v}^{\times \text{ }}}\mid f(x,y)\mid
_{v}^{s}\mid dxdy\mid  \notag \\
&=&\sum\limits_{m=1}^{\infty }q^{-(a+b)m\text{ }-d_{0\text{ }%
}ms}\int\limits_{^{O_{v}^{\times 2}\text{ }}}\mid f^{(m)}(x,y)\mid
_{v}^{s}\mid dxdy\mid ,  \label{I4}
\end{eqnarray}

with 
\begin{equation*}
f^{(m)}(x,y)=\pi ^{-md_{0}}f(\pi ^{am}x,\pi ^{bm}y)=\sum_{j=0}^{l_{f}}\pi
^{(d_{j}-d_{0})m}f_{j}(x,y).
\end{equation*}

The integral

\begin{equation*}
I(f^{(m)},s,v)=\int\limits_{O_{v}^{\times 2}\text{ }}\mid f^{(m)}(x,y)\mid
_{v}^{s}\mid dxdy\mid
\end{equation*}

is equal to 
\begin{equation}
I(s,f^{(m)},v,O_{v}^{\times 2})=\frac{U_{1}(q^{-s})}{1-q^{-1-s}}%
+\sum\limits_{\{\theta \in O_{v}^{\times }\mid f_{0}(1,\theta
^{a})=0\}}J(s,m,\theta ),  \label{I5}
\end{equation}
for $m$ big enough (cf. \ref{fm0a}). Thus from (\ref{I4}), and (\ref{I5}),
it follows that 
\begin{eqnarray}
Z(s,f,v,\Delta ) &=&\frac{U(q^{-s})}{1-q^{-1-s}}+  \notag \\
&&\sum\limits_{\{\theta \in O_{v}^{\times }\mid f_{0}(1,\theta
^{a})=0\}}\left( \sum\limits_{m=c(f)}^{\infty }q^{-(a+b)m\text{ }-d_{0\text{ 
}}ms}J(s,m,\theta )\right) ,  \label{I6}
\end{eqnarray}
where the constant $c(f)$ is defined in proposition \ref{prop3}. By using
the explicit formula for $J(s,m,\theta )$ given in proposition \ref{prop3}
for $m\geqslant c(f)$, and the algebraic identities given in subsection \ref
{saidentities},\ we have that

\begin{eqnarray}
\sum\limits_{m=c(f)}^{\infty }q^{-(a+b)m\text{ }-d_{0\text{ }%
}ms}J(s,m,\theta ) &=&\frac{U_{\theta }(q^{-s})}{\left( 1-q^{-1-s}\right) } 
\notag \\
&&+(1-q^{-1})^{2}\sum_{i=0}^{r-1}\frac{q^{-(s\mathcal{E}_{i+1}+1)}}{1-q^{-(s%
\mathcal{E}_{i+1}+1)}}S_{1}(s,i,c_{0})  \notag \\
&&-(1-q^{-1})^{2}\sum_{i=0}^{r-1}\frac{1}{1-q^{-(1+s\mathcal{E}_{i+1})}}%
S_{2}(s,i,c_{0})  \notag \\
&&+(1-q^{-1})^{2}\frac{q^{-(1+s\mathcal{E}_{r+1})}}{1-q^{-(1+s\mathcal{E}%
_{r+1})}}S_{1}(s,r,c_{0})  \notag \\
&&+\sum_{i=1}^{r}I_{i}(s)S_{2}(s,i-1,c_{0}),  \label{I7}
\end{eqnarray}
where $U_{\theta }(q^{-s})\in \mathbb{Q}\left[ q^{-s}\right] $, and the data 
$\mathcal{E}_{i}$, $\mathcal{D}_{i}$, $i=1,2,\cdots ,r$, depend on the
arithmetic Newton polygon $\ \ \Gamma _{f,\theta \text{.}_{\text{ }}}$Now by
remark \ref{remmark2}, the real parts of the poles of \ $S_{1}(s,i,c_{0})$, $%
S_{2}(s,i,c_{0})$, \ belong to the set \ $\mathcal{P}(\Gamma _{f,\theta })$.
Then the real parts of the poles of 
\begin{equation*}
\sum\limits_{\{\theta \in O_{v}^{\times }\mid f_{0}(1,\theta
^{a})=0\}}\left( \sum\limits_{m=c(f)}^{\infty }q^{-(a+b)m\text{ }-d_{0\text{ 
}}ms}J(s,m,\theta )\right)
\end{equation*}
\ belong to the set $\left\{ -\frac{a+b}{d_{0}}\right\} \bigcup
\bigcup\limits_{\{\theta \in O_{v}\mid f_{0}(1,\theta ^{a})=0\}}\mathcal{P}%
(\Gamma _{f,\theta })$, and from (\ref{I6}) follows that the real parts of
the poles of $Z(s,f,v,\Delta )$ belong to the set 
\begin{equation*}
\{-1\}\bigcup \{-\frac{a+b}{d_{0}}\}\bigcup \bigcup\limits_{\{\theta \in
O_{v}^{{}}\mid f_{0}(1,\theta ^{a})=0\}}\mathcal{P}(\Gamma _{f,\theta }).
\end{equation*}

\section{Main result}

In this section we prove the main result of this paper \ that gives an
explicit list \ for the possible poles of a local zeta function attached to
an arithmetically non degenerate polynomial in terms of the corresponding
arithmetic Newton polygon.

\subsection{Local zeta functions for arithmetically non degenerate curves}

Let $f(x,y)\in L_{v}[x,y]$ be a non-constant polynomial satisfying $f(0,0)=0$%
, and 
\begin{equation}
\mathbb{R}_{+}^{2}=\{(0,0)\}\bigcup \bigcup_{\gamma \subset \Gamma
^{geom}(f)}\Delta _{\gamma },
\end{equation}
\ a\ simplicial conical subdivision subordinated to $\Gamma ^{geom}(f)$. We
denote by $a_{\gamma }=$ $(a_{1}(\gamma ),a_{2}(\gamma ))$ a perpendicular
and primitive vector to facet $\gamma $ of $\Gamma ^{geom}(f)$, \ and \ by $%
<a,x>=d_{a}(\gamma )$ the equation of the corresponding supporting line. We
also define 
\begin{equation*}
\mathcal{P}(\Gamma ^{geom}(f)):=\{-\frac{a_{1}(\gamma )+a_{2}(\gamma )}{%
d_{a}(\gamma )}\mid \text{ }\gamma \text{ a facet of }\Gamma ^{geom}(f)\text{%
, with }d_{a}(\gamma )\neq 0\}.
\end{equation*}

The following is the main result of the present paper.

\begin{theorem}[Main Theorem]
\label{main1}Let $f(x,y)\in L_{v}\left[ x,y\right] $ be a non-constant
polynomial. If $\ f(x,y)$\ is arithmetically non degenerate with respect to
its arithmetic Newton polygon $\Gamma ^{A}(f)$, \ then the real parts of the
poles of \ $Z(s,f,v)$ belong to the set 
\begin{equation}
\{-1\}\bigcup \mathcal{P}(\Gamma ^{geom}(f))\bigcup \mathcal{P}(\Gamma
^{A}(f)).
\end{equation}
\end{theorem}

\begin{proof}
By taking a simplicial conical subdivision, the computation of \ $Z(s,f,v)$
is reduced to the computation of integrals of type $Z(s,f,v,O_{v}^{\times 2})
$, $\ $and $Z(s,f,v,\Delta _{\gamma })$, with\ $\gamma $ a proper face of $%
\Gamma ^{geom}(f)$ (see subsection \ref{subsection}).\ By lemma \ref{lemma2}
the real parts of the poles of \ $Z(s,f,v,O_{v}^{\times 2})$ belong to the
set $\{-1\}$. The computation \ of the integrals \ $Z(s,f,v,\Delta _{\gamma
})$, for \ $\gamma $ a proper face of $\Gamma ^{geom}(f)$, involves two
cases, according if the semi-quasihomogeneous \ polynomial \ $f(x,y),$ with
respect to $a_{\gamma }=$ $(a_{1}(\gamma ),a_{2}(\gamma ))$, is
geometrically non degenerate or not. If $\Delta _{\gamma }$ is a
one-dimensional cone, and $f_{\gamma }(x,y)$ does not have singularities on $%
\left( L_{v}^{\times }\right) ^{2}$, then the real parts of the poles of $%
Z(s,f,v,\Delta _{\gamma })$ belong to the set 
\begin{equation}
\{-1\}\bigcup \{-\frac{a+b}{d_{\gamma }}\}\subseteq \{-1\}\bigcup \mathcal{P}%
(\Gamma ^{geom}(f)),  \label{can1}
\end{equation}
where $ax+by=d_{\gamma }$, $d_{\gamma }\neq 0$, is the equation of the
supporting line of the facet $\gamma $ (cf. lemma \ref{lemma3}). If $\Delta
_{\gamma }$ is a two-dimensional cone, $f_{\gamma }(x,y)$ is a monomial and
then it does not have singularities on the torus $\left( L_{v}^{\times
}\right) ^{2}$, by lemma \ref{lemma3}, $Z(s,f,v,\Delta _{\gamma })$ \ is a
entire function.

If $\Delta _{\gamma }$ is a one-dimensional cone, and $f_{\gamma }(x,y)$ has
singularities on $\left( L_{v}^{\times }\right) ^{2}$, then $f(x,y)$ \ is a
semi-quasihomogeneous arithmetically non degenerate polynomial, and thus the
real parts of the poles of $Z(s,f,v,\Delta _{\gamma })$ belong to the set 
\begin{equation}
\{-1\}\bigcup \{-\frac{a+b}{d_{\gamma }}\}\bigcup \mathcal{P}(\Gamma
^{A}(f))\subseteq \{-1\}\bigcup \mathcal{P}(\Gamma ^{geom}(f))\bigcup 
\mathcal{P}(\Gamma ^{A}(f))\text{,}  \label{can2}
\end{equation}
(cf. theorem \ref{theorem1}).Therefore the real part of the poles of $%
Z(s,f,v)$ belong to the set $\{-1\}\bigcup \mathcal{P}(\Gamma
^{geom}(f))\bigcup \mathcal{P}(\Gamma ^{A}(f))$.
\end{proof}

\end{document}